\documentclass[11pt]{article}
\usepackage{amssymb}
\usepackage{amsmath}
\usepackage{amsthm}
\newcommand{\cA}{{\cal A}}
\newcommand{\cB}{{\cal B}}

\newcommand{\cH}{{\cal H}}

\newcommand{\cG}{{\cal G}}
\newcommand{\cO}{{\cal O}}
\newcommand{\cL}{{\cal L}}
\newcommand{\cM}{{\cal M}}
\newcommand{\cN}{{\cal N}}
\newcommand{\cF}{{\cal F}}
\newcommand{\cK}{{\cal K}}
\newcommand{\cP}{{\cal P}}
\newcommand{\cQ}{{\cal Q}}
\newcommand{\cR}{{\cal R}}
\newcommand{\cS}{{\cal S}}

\newcommand{\cU}{{\cal U}}
\newcommand{\cV}{{\cal V}}
\newcommand{\cW}{{\cal W}}
\newcommand{\cX}{{\cal X}}
\newcommand{\cY}{{\cal Y}}
\newcommand{\cZ}{{\cal Z}}

\newcommand{\NN}{{\mathbb N}}
\newcommand{\ZZ}{{\mathbb Z}}
\newcommand{\QQ}{{\mathbb Q}}

\newcommand{\gm}{\mathfrak{m}}    

\newcommand{\gp}{\mathfrak{p}}
\newcommand{\gq}{\mathfrak{q}}

\newcommand{\on}{\operatorname}
\newcommand{\Sh}{\on{Sh}}
\newcommand{\Fr}{\on{Fr}}
\newcommand{\EXT}{{{\cal E}xt}}
\newcommand{\Fl}{{\cal F}l}
\newcommand{\<}{\langle}
\renewcommand{\>}{\rangle}

\newcommand{\Qlb}{\mathbb{\bar Q}_\ell}

\newcommand{\A}{\mathbb{A}}
\newcommand{\Ql}{\mathbb{Q}_\ell}
\newcommand{\toup}[1]{\stackrel{#1}{\to}}
\newcommand{\hook}[1]{\stackrel{#1}{\hookrightarrow}}
\newcommand{\df}{\stackrel{\rm{def}}{=}}
\newcommand{\getsup}[1]{\stackrel{#1}{\gets}}

\newcommand{\IC}{\on{IC}}
\newcommand{\Hom}{\on{Hom}}
\newcommand{\Mod}{\on{Mod}}
\newcommand{\Ext}{\on{Ext}}

\newcommand{\Sym}{\on{Sym}}

\newcommand{\res}{\on{res}}
\newcommand{\Aut}{\on{Aut}}

\newcommand{\RG}{\on{R\Gamma}}

\newcommand{\diag}{\on{diag}}
\newcommand{\sym}{\on{sym}}

\newcommand{\triv}{\on{triv}}

\renewcommand{\part}{\on{part}}
\newcommand{\Spr}{{{\cal S}pr}}

\newcommand{\Pic}{\on{Pic}}

\newcommand{\Bun}{\on{Bun}}
\renewcommand{\Mod}{\on{Mod}}

\newcommand{\rk}{\on{rk}}
\newcommand{\Spec}{\on{Spec}}
\newcommand{\Specf}{\on{Spf}}

\newcommand{\supp}{\on{supp}}

\newcommand{\HOM}{{{\cal H}om}}
\newcommand{\END}{{{\cal E}nd}}

\newcommand{\GL}{\on{GL}}

\newcommand{\pr}{\on{pr}}
\newcommand{\id}{\on{id}}
\newcommand{\norm}{{norm}} 

\newcommand{\tr}{\on{tr}}
\newcommand{\graded}{{\on{gr}}}
\newcommand{\QED}{$\square$} 
\newcommand{\Fq}{\mathbb{F}_q}  
\newcommand{\Fp}{\mathbb{F}_p}  
\newcommand{\iso}{{\widetilde\to}}

\newcommand{\comp}{\circ}

\renewcommand{\H}{\on{H}}   
\newcommand{\R}{\on{R}\!}   

\newcommand{\D}{\on{D}}       
\newcommand{\wt}{\widetilde}
\newcommand{\select}[1]{{\it{#1}}}
\newcommand{\p}{\prime}

\newcommand{\suml}{\mathop{\sum}\limits}
\newcommand{\oplusl}{\mathop{\oplus}\limits}

\newcommand{\Gr}{{{\cG}r}}

\renewcommand{\div}{\on{div}}
\renewcommand{\P}{\mathbb{P}}

\newtheorem{Lm}{Lemma}
\newtheorem{Th}{Theorem}

\newtheorem{MLTh}{Main Local Theorem}

\newtheorem{MLThn}{Main Local Theorem$_n$}

\newtheorem{Pp}{Proposition}
\newtheorem{Cor}{Corolary}

\newtheorem{Slm}{Sublemma}

\theoremstyle{remark}
\newtheorem{Rem}{Remark}

\theoremstyle{definition}
\newtheorem{Def}{Definition}

\newenvironment{Prf}{\par\noindent {\it Proof }}{\QED}
\newcommand{\Step}[1]{\par\noindent{\bf Step {#1}}.}

\begin{document}
\author{Sergey Lysenko}
\title{Local geometrised Rankin-Selberg method for $\GL(n)$}
\date{}
\maketitle
\begin{abstract}
\noindent {\scshape Abstract} \hskip 0.8 em 
Following Laumon \cite{L1}, to a nonramified $\ell$-adic local system $E$
of rank $n$ on a curve $X$ one associates a complex of $\ell$-adic sheaves
$_n\cK_E$ on the moduli stack of rank $n$ vector bundles on $X$ with
a section, which is cuspidal and satisfies Hecke property for $E$.
This is a geometric counterpart of the well-known construction due 
to Shalika \cite{Shal} and Piatetski-Shapiro \cite{P-Sh}.  
We express the cohomology of the tensor product 
$_n\cK_{E_1}\otimes {_n\cK_{E_2}}$ in terms of cohomology of 
the symmetric powers of $X$. This may be considered as a geometric
interpretation of the local part of the classical Rankin-Selberg 
method for $\GL(n)$ in the framework of the geometric Langlands program.
\end{abstract}

\subsection{General introduction}

This is the first in a series of two papers, where we propose a geometric
version of the classical Rankin-Selberg method for computation of the scalar
product of two cuspidal automorphic forms on $\GL(n)$ over a function field. 
This geometrization fits in the framework of the geometric Langlands program 
initiated by V. Drinfeld, A. Beilinson and G. Laumon.

 Let $X$ be a smooth, projective, geometrically connected
curve over $\Fq$. Let $\ell$ be a prime invertible in $\Fq$. 
According to the Langlands correspondence for $\GL(n)$ over 
function fields (proved by L. Lafforgue), 
to any smooth geometrically irreducible $\Qlb$-sheaf $E$ of rank $n$
on $X$ is associated a (unique up to a multiple) cuspidal 
automorphic form $\varphi_E: \Bun_n(\Fq)\to \Qlb$, which
is a Hecke eigenvector with respect to $E$. The function $\varphi_E$ 
is defined on the set $\Bun_n(\Fq)$ of isomorphism classes 
of rank $n$ vector bundles on $X$. 

 The classical method of Rankin and Selberg for $\GL(n)$ may be divided into
two parts: local and global. The global result calculates 
for any integer $d$ the scalar product of two (appropriately normalized)
automorphic forms
\begin{equation}
\label{scalar_square}
\sum_{L\in\Bun^d_n(\Fq)} \frac{1}{\#\Aut L}\; 
\varphi_{E_1^*}(L)\varphi_{E_2}(L),
\end{equation}
where $\Bun^d_n(\Fq)$ is the set of isomorphism classes of vector 
bundles $L$ on $X$ of rank $n$ and degree $d$, and $\#\Aut L$ 
stands for the number of elements in $\Aut L$. More precisely, 
this scalar product vanishes if and only if $E_1$ and $E_2$ are non 
isomorphic. In the case $E_1\iso E_2\iso E$ the answer is expressed 
in terms of the action of the geometric Frobenius endomorphism on 
$\H^1(X\otimes\bar\Fq,\,\END E)$.

The computation of (\ref{scalar_square}) is based on the equality of 
formal series
\begin{equation}
\label{equality_formal_series}
\sum_{d\ge 0}\;\; \sum_{(\Omega^{n-1}\hook{}L)\in
\,{_n\cM_d}(\Fq)}\;\;\frac{1}{\#\Aut(\Omega^{n-1}\hook{}L)}
\varphi_{E_1^*}(L)\varphi_{E_2}(L)t^d=
L(E_1^*\otimes E_2, \, q^{-1}t)
\end{equation}
Here $_n\cM_d(\Fq)$ is the set of isomorphism classes of pairs
$(\Omega^{n-1}\hook{}L)$, where $L$ is a vector bundle on $X$ 
of rank $n$ and degree $d+n(n-1)(g-1)$, and $\Omega$ is the
canonical invertible sheaf on $X$ ($\Omega^{n-1}$ is embedded in $L$ 
as a subsheaf, i.e., the quotient is allowed to have torsion). 
We have denoted by $L(E_1^*\otimes E_2,t)$ 
the L-function attached to the local system $E_1^*\otimes E_2$ on $X$. 

 Recall that the existence of the automorphic form $\varphi_E$ is a descent
problem (cf. \cite{L1}). Using an explicit construction due 
to Shalika \cite{Shal} and Piatetski-Shapiro \cite{P-Sh}, one
associates to a smooth $\Qlb$-sheaf $E$ of rank $n$ on $X$  
a function $\tilde\varphi_E: {_n\cM_d}(\Fq)\to\Qlb$, which is cuspidal
and satisfies Hecke property with respect to $E$. The 
Langlands conjecture predicts that when $E$ is geometrically
irreducible, $\tilde\varphi_E$ is constant along the fibres of 
the projection $_n\cM_d(\Fq)\to \Bun_n^{d+n(n-1)(g-1)}(\Fq)$,
that is, $\tilde\varphi_E$ is the pull-back of a function 
$\varphi_E$ on $\Bun_n(\Fq)$. So, (\ref{equality_formal_series}) 
is a statement independent of the Langlands conjecture. 
In fact, (\ref{equality_formal_series}) is 
of local nature: it is true for any local systems $E_1$ and $E_2$ 
of rank $n$ on $X$ after replacing $\varphi_E$ by $\tilde\varphi_E$. 

 Main result of this paper is a strengthened geometric version of
the equality
\begin{equation}
\label{equality_coefficients}
\sum_{(\Omega^{n-1}\hook{}L)\in
\,{_n\cM_d}(\Fq)}\;\;\frac{1}{\#\Aut(\Omega^{n-1}\hook{}L)}
\tilde\varphi_{E_1^*}(L)\tilde\varphi_{E_2}(L)=
q^{-d}\!\!\!\!\!\sum_{D\in X^{(d)}(\Fq)} 
\tr(\Fr, (E_1^*\otimes E_2)^{(d)}_D)
\end{equation}
of coefficients in (\ref{equality_formal_series}) for each $d\ge 0$.
Here $X^{(d)}$ is the $d$-th symmetric power of $X$,  $(E_1^*\otimes
E_2)^{(d)}$ is a constructable $\Qlb$-sheaf on $X^{(d)}$ (cf.
Sect.~\ref{Sect_Springer}), and $\Fr$ is the geometric Frobenius
endomorphism.  

Let $_n\cM_d$ denote the moduli stack of pairs 
$(\Omega^{n-1}\hook{s}L)$, where $L$ is a vector bundle 
of rank $n$ and degree $d+n(n-1)(g-1)$ on $X$, and $s$ is an inclusion
of $\cO_X$-modules. Following Drinfeld \cite{Dr} (n=2) and Deligne (n=1), 
Laumon \cite{L1} has defined a complex of $\Qlb$-sheaves $_n\cK^d_E$
on $_n\cM_d$, which is
a geometric counterpart of $\tilde\varphi_E$\footnote{We normalize
$_n\cK^d_E$ as in Remark~\ref{Rem_functors_well-defined} (cf. Sect.~2.1).
This also gives a normalization of $\tilde\varphi_E$ as the function `trace
of Frobenius' of $_n\cK^d_E$.}. The geometric
Langlands conjecture predicts that when $E$ is a smooth geometrically
irreducible $\Qlb$-sheaf of rank $n$ on $X$, $_n\cK^d_E$ descends with
respect to the projection $_n\cM_d\to\Bun_n$, where $\Bun_n$ is the moduli
stack of rank $n$ vector bundles on $X$. 

 We establish for any smooth $\Qlb$-sheaves $E_1,E_2$ of rank $n$ on $X$ and
any $d\ge 0$ a canonical isomorphism 
$$
\RG_c(_n\cM_d\, ,\, {_n\cK^d_{E_1^*}}\otimes{_n\cK^d_{E_2}})\iso
\RG( X^{(d)},\, (E_1^*\otimes E_2)^{(d)})(d)[2d],
$$
which is a geometric version of (\ref{equality_coefficients}). In fact, a
more general statement is proved.

\medskip
\centerline{\bf Acknowledgements}

\smallskip

I am deeply grateful to V. Drinfeld, who has posed the problem we 
consider and explained to me its place in the geometric Langlands 
program. In this paper we attempt to realize certain of his ideas. 
I am grateful to D. Gaitsgory who has suggested to the author an idea of the
second proof of Theorem~\ref{Th_B} (cf. Sect. 6.6).
I also wish to thank G.~Laumon for constant support.

\subsection{Conventions and notation}
0.2.1. Fix an algebraically closed ground field $k$ of characteristic
$p>0$, a prime $\ell\ne p$ and an algebraic closure $\Qlb$ of $\Ql$. 
All the schemes and stacks we use will be defined over $k$. 
Throughout the paper, $X$ will denote a fixed smooth projective
connected curve of genus $g\ge 1$ (over $k$). 

 We will work with algebraic stacks in smooth topology and with
(perverse) $\Qlb$-sheaves on them.
If $\cX$ is an algebraic stack  locally of finite type then the 
notion of a (perverse) $\Qlb$-sheaf on $\cX$ localizes in the 
smooth topology, and hence makes perfect sense. However, the
corresponding derived category is problematic. We adopt the point
of view that an appropriate formalism exists (it is partially 
established in \cite{LMb}). Let $f:\cX\to\cY$ be a 
morphism of algebraic stacks. The functors $f^*, f_*,
f_!$ will be understood in the derived category sense. 

 We say that $f$ is a
\select{generalized affine fibration of rank $m$} in the following cases.
First, if locally in smooth topology on $\cY$  there exists a
homomorphism
$L\to L'$ of locally free coherent sheaves  on $\cY$ and an $L'$-torsor
$\cY'\to\cY$ such that $f$ is identified  with $\cY'/L\to \cY$, the quotient
being taken in stack sense, and 
$\rk L'-\rk L=m$. Second, if the map $f$ can be written as the composition
of generalized affine fibrations of first type of ranks 
$m_1,\ldots,m_k$ with $\sum m_i=m$. We essentially use the fact 
that for a generalized affine fibration
$f$ of rank $m$ one has $f_!\,\Qlb\;\iso\;\Qlb(-m)[-2m]$. 

  We fix a nontrivial
additive  character $\psi:\Fp\to \Qlb^*$ and denote  by $\cL_{\psi}$ the
Artin-Schreier sheaf on $\A^1_k$ associated to $\psi$ (SGA$4\frac{1}{2}$, 
[Sommes trig.], 1.7). Fix also a square root of $p$ in $\Qlb$ and
define using it the sheaf $\Qlb(\frac{1}{2})$ over $\Spec \Fp$ and, 
hence, over $\Spec k$.

\medskip
\noindent
0.2.2. When we say that a stack $\cY$ classifies \select{something}, 
it should always be clear what an $S$-family of \select{something} 
is for any $k$-scheme $S$, i.e., what is the groupo\"\i d $\Hom(S,\cY)$ 
and what are the functors $\Hom(S_2,\cY)\to\Hom(S_1,\cY)$
for each morphism $S_1\to S_2$.

 For example, if $\cY$ is the stack that classifies pairs $M_1\hook{} M_2$
with $M_1$ (resp. $M_2$) being a coherent sheaf on $X$ of generic rank 
$i_1$ and of degree $d_1$ (resp., of generic rank $i_2$ and of degree 
$d_2$) then $\Hom(S,\cY)$ is the groupo\"\i d whose objects are 
inclusions $M_1\hook{} M_2$ of coherent sheaves on 
$S\times X$ that are $S$-flat and such that the quotient $M_2/M_1$ 
is also $S$-flat, and for any point $s\in S$ the conditions on 
the generic rank and on the degree of $M_i\mid_{s\times X}$ $(i=1,2)$ 
hold. Morphisms from an object $M_1\hook{} M_2$ to
an object $M'_1\hook{} M'_2$ are by definition the isomorphisms 
$M_1\iso M'_1$ and $M_2\iso M'_2$ making the natural diagram commutative. 

We denote by $\Sh_i$ the moduli stack of coherent sheaves on $X$ 
of generic rank $i$. This is an algebraic stack locally of finite type. 
Its connected components are numbered by $d\in \ZZ$:
the component $\Sh^d_i$ classifies coherent sheaves of rank $i$ and of degree
$d$ on $X$. The stack $\Sh^d_0$ is, in fact, of finite type. 

 By $\Pic X\subset\Sh_1$ we denote the open substack classifying
invertible $\cO_X$-modules. This is the Picard stack of $X$. Its
connected component $\Pic^d X$ classifies line bundles of degree $d$ on $X$. 

 Denote by $^{\le n}\Sh^d_0\subset \Sh^d_0$
the open substack given by the property: for a scheme $S$ an object $F$ of
$\Hom(S, \Sh^d_0)$ lies in $\Hom(S, {^{\le n}\Sh^d_0})$ if the geometric 
fibre of $F$ at any point of $X\times S$ is of dimension at most $n$. 
We write $X^{(d)}$ for the $d$-th symmetric power of $X$.
By $\div:\Sh^d_0\to X^{(d)}$ is denoted
the morphism norm (cf. \cite{Gr},6). If $D_1,\dots,D_s$ are effective divisors
on $X$ then it sends the $\cO_X$-module $\cO_{D_1+\ldots+D_s}\oplus
\cO_{D_2+\ldots+D_s}\oplus\ldots\oplus\cO_{D_s}$ to $D_1+2D_2+\ldots+sD_s$.

\medskip
\noindent
0.2.3.  Fix the maximal
torus of diagonal matrices in $\GL(n)$ and the Borel subgroup 
of upper-triangular matrices. Then the set of weights of $\GL(n)$ is
identified with $\ZZ^n$. The fundamental weights are given by 
$\omega_i=(1,\ldots,1,0,\ldots,0)\in\ZZ^n$, where 1 occurs $i$ times
($i=1,\ldots,n$). 

Define the following semigroups 
$\Lambda^+_n\subset\Lambda_n\subset \Lambda^p_n$ consisting of weights.  
Let $\Lambda_n=\ZZ_+^n$ and
$\Lambda^p_n=\{\lambda\in\ZZ^n \mid \lambda_1+\ldots+\lambda_i\ge 0\;\,
\mbox{for all}\;\, i\}$. The superscript $p$ should designate
that $\Lambda^p_n$ contains the $\ZZ_+$-span of positive roots. 
Set also $\Lambda_n^+=\{\lambda=
(\lambda_1\ge\ldots\ge\lambda_n\ge 0) \mid
\lambda_i\in\ZZ\}$. Similarly, we let
$\Lambda_n^-=\{\lambda=(0\le\lambda_1\le\ldots\le
\lambda_n) \mid \lambda_i\in\ZZ\}$. 

 For $d\ge 0$ we also introduce $\Lambda_{n,d}\subset\Lambda_n$, 
$\Lambda_{n,d}^+\subset\Lambda_n^+$ and so on, where the subscript $d$
means that we impose the condition $\sum \lambda_i=d$. 
The half sum of positive roots is denoted by $\rho$. 

 For a weight $\lambda$ of $\GL(n)$ we introduce the schemes $X^{\lambda}_+,
X^{\lambda}_-, X^{\lambda}$ and $X^{\lambda}_p$ that should be thought of as 
the moduli schemes of $\Lambda_n^+$ (resp., of $\Lambda_n^-$, $\Lambda_n$, $\Lambda_n^p$)
\  -valued divisors on $X$ of degree $\lambda$. The precise definition is
as follows. 

 Set $X^{\lambda}_p=\prod_{i=1}^n X^{(\lambda_1+\ldots+\lambda_i)}$.
A point of $X^{\lambda}_p$ is a collection of (not necessarily effective) 
divisors $(D_1,\ldots,D_n)$ on $X$ with $D_1+\ldots+D_i\in 
X^{(\lambda_1+\ldots+\lambda_i)}$. Let $X^{\lambda}\hook{} X^{\lambda}_p$ 
be the closed subscheme given by $D_i\ge 0$ for all $i$. Let $X^{\lambda}_+$ (resp.,
$X^{\lambda}_-$) be the closed subscheme of $X^{\lambda}$ given by
$D_1\ge\ldots\ge D_n$ (resp, $D_1\le\ldots\le D_n$).

\smallskip

 Given a closed point $(D_i)$ of $X^{\lambda}_i$ with $D_i=\sum_x d_{i,x} x$,
we associate to it a divisor on $X$ with values in $\Lambda^p_n$. The value
of this divisor at $x$ is the weight $(d_{1,x},\ldots,d_{n,x})$. In the same way
a closed point of $X^{\lambda}$ (resp., $X^{\lambda}_+$,$X^{\lambda}_-$) can be
viewed as a $\Lambda_n$ (resp., $\Lambda^+_n$,$\Lambda^-_n$) -valued divisor
on $X$.

\medskip
\noindent
0.2.4. For $\lambda\in \Lambda_{n,d}^+$ define the polynomial functor
$V^{\lambda}$ of a $\Qlb$-vector space $V$ as follows. Let
$\lambda=(\lambda_1,\ldots,\lambda_{n'},0,\ldots,0)$ with
$\lambda_{n'}>0$. Denote by $U^{\lambda}$ the irreducible
representation of $S_d$ (over $\QQ$) associated to $\lambda$. So, for
example, if $\lambda=(d,0,\ldots,0)$ then $U^{\lambda}\iso\QQ$ is
trivial, and if $\lambda=(1,\ldots,1)$ then $U^{\lambda}$ is the
signature representation. Set 
$$
V^{\lambda}=(V^{\otimes d}\otimes_{\QQ} U^{\lambda})^{S_d},
$$
where it is understood that $S_d$ acts by permutations on $V^{\otimes
d}$ and diagonally on the tensor product. If $m=\dim V< n'$ then
$V^{\lambda}=0$, otherwise $V^{\lambda}$ is the irreducible
representation of $\GL(V)$ of the highest weight
$(\lambda_1,\ldots,\lambda_{n'},0,\ldots,0)\in\Lambda_{m,d}^+$.

\section{Laumon's perverse sheaf $\cL^d_E$}
\label{Sect_Springer}
Let $E$ be a smooth $\Qlb$-sheaf on $X$. Recall the definition 
of Laumon's perverse sheaf $\cL^d_E$ on $\Sh^d_0$ 
associated to $E$ (\cite{L1}). 
Denote by $\sym: X^d\to X^{(d)}$ the natural map and  consider the smooth
$\Qlb$-sheaf $E^{\boxtimes\, d}$ on $X^d$. Notice that 
$\sym_!(E^{\boxtimes\, d})[d]$ is a perverse sheaf.  Set
$$
E^{(d)}=(\sym_!(E^{\boxtimes\, d}))^{S_d}
$$ 
Since $E^{(d)}$ is a direct 
summand of $\sym_!(E^{\boxtimes\, d})$, \ $E^{(d)}[d]$
is also a perverse sheaf. 
 
Denote by $\Fl^{1,\ldots,1}$ (1 occurs $d$ times) the stack of complete flags
$(F_1\subset\ldots\subset F_d)$, where $F_i$ is a coherent torsion sheaf on $X$ of
length $i$. The morphism $\gp: \Fl^{1,\ldots,1}\to\Sh^d_0$ that sends 
$(F_1\subset\ldots\subset F_d)$ to $F_d$ is
representable and proper. The morphism
$\gq:\Fl^{1,\ldots1}\to\Sh^1_0\times\ldots\times\Sh^1_0$ that sends
$(F_1\subset\ldots\subset F_d)$ to $(F_1, F_2/F_1,\ldots, F_d/F_{d-1})$ is a
generalized affine fibration. This, in particular, implies that $\Fl^{1,\ldots,1}$ is
smooth. 

 Springer's sheaf $\Spr^d_E$ on $\Sh^d_0$ is defined as
$$
\Spr^d_E=\gp_! \gq^* (\div^{\times d})^* (E^{\boxtimes\, d})
$$
Since $\gp$ is small, $\Spr^d_E$ is a perverse sheaf that 
coincides with the Goresky-MacPherson extension of its restriction to any 
nonempty open substack of $\Sh^d_0$. It also carries a natural $S_d$-action
(cf. Theorem 3.3.1 of \cite{L1}). 
Set 
$$
\cL^d_E=\Hom_{S_d}(\triv, \Spr^d_E),
$$
where $\triv$
denotes the trivial representation of the symmetric group $S_d$.
Again, $\cL^d_E$ is a direct summand of $\Spr^d_E$, so $\cL^d_E$ is
perverse and coincides with the
Goresky-MacPherson extension of its restriction to any nonempty open substack 
of $\Sh^d_0$. We have a smooth morphism $X^{(d)}\to\Sh^d_0$ that sends a divisor
$D$ to $\cO_D$, and the pull-back of $\cL^d_E$ under this map 
is identified with $E^{(d)}$.

\section{Main results}
2.1 \  Fix $n>0, d\ge 0$. Let $\Omega$ be the canonical invertible
sheaf on $X$. Denote by
$_n\cQ_d$ the stack that classifies collections 
\begin{equation}
\label{collect_1}
(0=L_0\subset L_1\subset \ldots\subset L_n\subset L, \;\;\; (s_i)),
\end{equation}
where $L_n\subset L$ is a modification of rank $n$ vector bundles on $X$
with $\deg(L/L_n)=d$, $(L_i)$ is a complete flag of subbundles on $L_n$,
and $s_i:\Omega^{n-i}\iso L_i/L_{i-1}$ is an isomorphism 
$(i=1,\ldots, n)$. We have a map $\mu: {_n\cQ_d}\to \A^1_k$ which at the
level of $k$-points sends the above collection to the sum 
of $n-1$ classes in
$$
k\iso \Ext^1(\Omega^{n-i-1}, \Omega^{n-i})\iso \Ext^1(L_{i+1}/L_i, 
L_i/L_{i-1})
$$
that correspond to the succesive extensions 
$
0\to  L_i/L_{i-1}\to L_{i+1}/L_{i-1}\to L_{i+1}/L_i\to 0
$.

Let $\beta:{_n\cQ_d}\to {^{\le n}\Sh^d_0}$ be the map that sends 
(\ref{collect_1}) to
$L/L_n$. It is of finite type and smooth of relative dimension
$
b=b(n,d)=nd+(1-g)\suml_{i=1}^{n-1} i^2
$. Therefore,  $_n\cQ_d$ is smooth and of finite type.  
So, if $E$ is a smooth
$\Qlb$-sheaf on $X$ then on $_n\cQ_d$ we have a perverse $\Qlb$-sheaf
$$
_n\cF^d_{E,\psi}=\beta^*\cL^d_E\otimes
\mu^*\cL_{\psi}[b](\frac{b}{2})
$$

 Let $\pi_0: {_n\cQ_d}\to X^{(d)}$ be the map 
that sends (\ref{collect_1}) to the divisor $D\in X^{(d)}$ for which
the inclusion of invertible sheaves $\wedge^n L_n\hook{} \wedge^n L$
induces an isomorphism $\wedge^n L_n(D)\iso \wedge^n L$. We also
have a map $X^{(d)}\to \Pic^d X$ that sends a divisor $D$ to $\cO_X(D)$.

Let $_n\cM_d$ be the stack 
classifying pairs
$(\Omega^{n-1}\hook{} L)$, where $L$ is an
$n$-bundle on $X$
with 
$$
\deg L- \deg (\Omega^{(n-1)+(n-2)+\ldots+(n-n)})=d
$$ 
The forgetful map $\zeta: {_n\cQ_d}\to {_n\cM_d}$ is
representable, and the following diagram commutes:
$$
\begin{array}{ccc}
_n\cQ_d & \toup{\pi_0} & X^{(d)}\\
\downarrow\lefteqn{\scriptstyle \zeta} &&\downarrow\\
_n\cM_d & \toup{\theta} & \Pic^d X,
\end{array}
$$
where $\theta$ is the map that sends $(\Omega^{n-1}\hook{} L)$ to 
$\det L\otimes\Omega^{(1-n)+(2-n)+\ldots+(n-n)}$. Denote by 
$$
\pi: {_n\cQ_d}\times_{_n\cM_d}{_n\cQ_d}\to X^{(d)}\times_{\Pic^d X}
X^{(d)}
$$
the morphism $\pi_0\times\pi_0$. Since $X^{(d)}\to\Pic^d X$ is
representable and separated, the diagonal map $i: X^{(d)}\to
X^{(d)}\times_{\Pic^d X} X^{(d)}$ is a closed immersion. Our main result
is the next theorem.

\begin{MLTh}
For any smooth $\Qlb$-sheaves $E,E'$ on $X$ of ranks $m,m'$ respectively
with $\min\{m,m'\}\le n$ there exists a canonical isomorphism
$$
\pi_!({_n\cF^d_{E,\psi}}\boxtimes \, {_n\cF^d_{E',\psi^{-1}}})\,\iso\, 
i_*(E\otimes E')^{(d)}(d)[2d]
$$
in the derived category on $X^{(d)}\times_{\Pic^d X} X^{(d)}$. 
\end{MLTh}
\begin{Rem}
\label{Rem_functors_well-defined}
i) The stack $_n\cQ_d\times_{_n\cM_d}{_n\cQ_d}$ is of finite type, though
$_n\cM_d$ is not, so that $\pi$ is of finite type but not representable. 

\smallskip\noindent
ii) Define the complex $_n\cK^d_E$ on $_n\cM_d$ as
$_n\cK^d_E=\zeta_!(_n\cF^d_{E,\psi})$. The geometric Langlands conjecture
claims that if $E$ is a smooth irreducible $\Qlb$-sheaf of rank $n$ on $X$
then for each $d\ge 0$ the complex $_n\cK^d_E$ descends with respect to
the projection $_n\cM_d\to \Bun_n$. 
\end{Rem}

\medskip\smallskip
\noindent
2.2 \ Actually we prove a more general statement.
Recall that for $\Qlb$-vector spaces
$E, E'$  of dimensions $m, m'$ respectively we have 
$$
\Sym^d(E\otimes
E')=\oplus_{\lambda\in\Lambda_{r,d}^+}\; E^{\lambda}\otimes
(E')^{\lambda},
$$ 
where $r=\min\{m,m'\}$. To formulate the version of Main Local Theorem we
actually prove, we globalize the above equality as follows.   

 For $\lambda\in\Lambda_{n,d}^+$ and a smooth $\Qlb$-sheaf $E$ on $X$ we 
define a constructable $\Qlb$-sheaf $E^{\lambda}_+$ on $X^{\lambda}_+$ 
(cf. Sect. 3.1), which is a
global analog of the corresponding polynomial functor. The fibre of
$E^{\lambda}_+$ at $D=\sum_x \lambda_x x$ is the tensor product over
closed points of $X$
$$
\otimes_{x\in X}\; (E_x)^{\lambda_x},
$$
where $E_x$ denotes the fibre of $E$ at $x$. For example, for
$\lambda=(d,0,\ldots,0)$ we have $X^{\lambda}_+=X^{(d)}$ and $E^{\lambda}_+=
E^{(d)}$. Another example, for $\lambda=\omega_i$ we obtain $X^{\lambda}_+=X$
and $E^{\lambda}_+=\wedge^i E$. 
 
 Denote by $\pi^{\lambda}:
X^{\lambda}_+\to X^{(d)}$ the map that sends $(D_1\ge\ldots\ge D_n\ge 0)\in
X^{\lambda}_+$ to $\sum D_i$. 

\begin{Lm} 
\label{Lm_filtration}
For any smooth $\Qlb$-sheaves $E, E'$ on
$X$ of ranks $m,m'$ resp. there is a canonical filtration 
$$
0={^{\le 0}}(E\otimes E')^{(d)}\subset{^{\le 1}}(E\otimes E')^{(d)}\subset\ldots
$$ 
on $(E\otimes E')^{(d)}$ by constructable subsheaves with
the following property. First, if $\min\{m,m'\}\le n$ then 
$^{\le n}(E\otimes E')^{(d)}=(E\otimes E')^{(d)}$. Secondly,
there is a canonical refinement of this filtration such that
$$
\graded\; {^{\le n}(E\otimes E')^{(d)}}\,\iso\,
\oplus_{\lambda\in\Lambda_{n,d}^+}
\; \pi^{\lambda}_* (E^{\lambda}_+\otimes
E^{\p\,\lambda}_{\, +})
$$
for each $n$. 
\end{Lm}

\begin{MLThn} 
For any smooth $\Qlb$-sheaves $E,E'$ on $X$ there exists a canonical
isomorphism
$$
\pi_!({_n\cF^d_{E,\psi}}\boxtimes \, {_n\cF^d_{E',\psi^{-1}}})\,\iso\, 
i_*\, {^{\le n}(E\otimes E')^{(d)}}(d)[2d]
$$
in the derived category on $X^{(d)}\times_{\Pic^d X} X^{(d)}$. 
\end{MLThn}

\medskip
\noindent 
2.3 \  The proof consists of the following steps. 
Let us
denote by $_n\cX_d$ the stack classifying collections $(L, (t_i))$,
where $L$ is a vector bundle on $X$ of rank $n$, 
$$
t_i:\Omega^{(n-1)+(n-2)+
\ldots+(n-i)}\hook{} \wedge^i L
$$
is an inclusion of $\cO_X$-modules
$(i=1,\ldots, n)$, and $\deg L- \deg (\Omega^{(n-1)+(n-2)+\ldots+(n-n)})
=d$. 

Given an object of $_n\cQ_d$, we get the morphisms 
$$
t_i:\Omega^{(n-1)+\ldots+(n-i)}\iso \wedge^i L_i\hook{} \wedge^i L
$$
This defines a map $\varphi: {_n\cQ_d}\to{_n\cX_d}$. Notice that 
$\zeta: {_n\cQ_d}\to {_n\cM_d}$ factors as
$_n\cQ_d\toup{\varphi}{_n\cX_d}\to{_n\cM_d}$, 
where the second arrow is the forgetful map. 
Since $_n\cX_d\to{_n\cM_d}$ is representable and separated,
the natural map
$
{_n\cQ_d}\times_{_n\cX_d}{_n\cQ_d}\to
{_n\cQ_d}\times_{_n\cM_d}{_n\cQ_d}
$
is a closed immersion. Let
$$
\pi':{_n\cQ_d}\times_{_n\cX_d}{_n\cQ_d}\to X^{(d)}\times_{\Pic^d X}
X^{(d)}
$$  
be the restriction of $\pi$ to
${_n\cQ_d}\times_{_n\cX_d}{_n\cQ_d}$. 
The first step is to establish the following result.

\begin{Th}
\label{Th_A}
For any smooth $\Qlb$-sheaves $E,E'$ on $X$ the natural map
$$
\pi_!({_n\cF^d_{E,\psi}}\boxtimes \, {_n\cF^d_{E',\psi^{-1}}})\to
\pi'_!({_n\cF^d_{E,\psi}}\boxtimes \, {_n\cF^d_{E',\psi^{-1}}})
$$
is an isomorphism.
\end{Th}

 Our proof of Theorem~\ref{Th_A} will be based on
Proposition~\ref{Pp_Cor_CS}, which is a corolary of
the geometric Casselman-Shalika formula for $\GL(n)$ (cf. \cite{N1,FGV,N2}). 
We present it in Sect.~\ref{section_around_CS} written 
independently of the rest of the paper. 

\medskip

 The second step is as follows. Let $\phi:
{_n\cX_d}\to X^{(d)}$ be the map that sends $(L, (t_i))$ to the divisor $D
\in X^{(d)}$ such that $t_n$ induces an isomorphism
$$
\Omega^{(n-1)+\ldots+(n-n)}(D)\iso
\wedge^n L,
$$ 
so that $\phi\comp\varphi=\pi_0$. 
We will write $f: {_n\cQ_d}\times_{_n\cX_d}{_n\cQ_d} \to X^{(d)}$
for
the composition $_n\cQ_d\times_{_n\cX_d}{_n\cQ_d} 
\to {_n\cX_d}\toup{\phi} X^{(d)}$,
where the first map is the natural projection. The morphism $f$ is of finite
type but not representable. Since the diagram
$$
\begin{array}{ccc}
{_n\cQ_d}\times_{_n\cX_d}{_n\cQ_d} & \hook{} &
{_n\cQ_d}\times_{_n\cM_d}{_n\cQ_d}\\
\downarrow\lefteqn{\scriptstyle f} &&
\downarrow\lefteqn{\scriptstyle\pi}\\ X^{(d)} & \toup{i} &
X^{(d)}\times_{\Pic^d X} X^{(d)}
\end{array}
$$
commutes, Main Local Theorem is just a combination of Theorem~\ref{Th_A}
with the following result.

\begin{Th}
\label{Th_B}
For any smooth $\Qlb$-sheaves $E,E'$ on $X$ there is a canonical isomorphism
\begin{equation}
\label{iso_1}
f_!({_n\cF^d_{E,\psi}}\boxtimes \, {_n\cF^d_{E',\psi^{-1}}})\,\iso\,\,
{^{\le n}(E\otimes E')^{(d)}}(d)[2d]
\end{equation}
\end{Th}

\smallskip

We present two different proofs of Theorem~\ref{Th_B}. In the first proof,
which occupies Sect. 6.1 through 6.5, we derive Theorem~\ref{Th_B}
from the following result.

\begin{Th}
\label{Th_C}
Let ${^{\le n}\div}: {^{\le n}\Sh^d_0}\to X^{(d)}$ denote the restriction of
$\div:\Sh^d_0\to X^{(d)}$. For any smooth $\Qlb$-sheaves $E, E'$ on $X$ 
the complex $(^{\le n}\div)_!(\cL^d_E\otimes\cL^d_{E'})$ is placed in degrees 
$\le -2d$, and for the highest cohomology sheaf of this complex we have
canonically
$$
\R^{-2d}(^{\le n}\div)_! (\cL^d_E\otimes\cL^d_{E'})(-d)\,\iso \;{^{\le
n}(E\otimes E')^{(d)}}
$$
\end{Th}

\smallskip

 In Sect. 6.6 we present an alternative proof of Theorem~\ref{Th_B}.
The idea of this proof was communicated to the author by D. Gaitsgory.
This proof requires the additional assumption: $\min\{\rk E,\rk E'\}\le n$. 
The reader interested in the proof of Mail Local Theorem under this assumption
may skip Sect. 6.1 through~6.5.

\section{Around the geometric Casselman-Shalika formula for $\GL(n)$}
\label{section_around_CS}

3.1 \ The purpose of Sect.~\ref{section_around_CS} is to present
Proposition~\ref{Pp_Cor_CS}, which is a corolary
of the geometric Casselman-Shalika formulae for $\GL(n)$ (cf.
\cite{N1,FGV,N2}).  To formulate it we introduce some
notation.

 Fix $\lambda\in\Lambda_{n,d}^-$. Recall that
$X^{\lambda}_-$ is the scheme of collections $(D_1,\ldots,
D_n)$, where $D_i$ is an effective divisor on $X$ of degree $\lambda_i$
with $D_1\le\ldots\le D_n$. Let
$$
i_{\lambda}: {X^{\lambda}_-}\to {^{\le n}\Sh^d_0}
$$ 
be the map that sends
$(D_1,\ldots,D_n)$ to 
$$
\Omega^{n-1}(D_1)/\Omega^{n-1}\oplus
\Omega^{n-2}(D_2)/\Omega^{n-2}\oplus\ldots\oplus \cO(D_n)/\cO
$$ 
According to (\cite{L1}, Theorem 3.3.8),
if $E$ is a smooth $\Qlb$-sheaf on $X$ then the complex
$i_{\lambda}^*\cL^d_E$ is placed in degrees $\le 2a(\lambda)$ with respect
to the usual t-structure, where
$$
a(\lambda)\df<\lambda, (n-1, n-2, \ldots,0)>
$$
Moreover, if $m\in\NN$ is such that $\lambda=
(0,\ldots,0,\lambda_{n-m+1},\ldots,\lambda_n)$ with $\lambda_{n-m+1}>0$
then $2a(\lambda)$-th cohomology sheaf of $i_{\lambda}^*\cL^d_E$ vanishes
if and only if $\rk E<m$. 

 For a weight $\lambda=(\lambda_1,\ldots,\lambda_n)$ set $\lambda^t=
(\lambda_n,\ldots,\lambda_1)$. Denote also by $t: X^{\lambda}_-\iso
X^{\lambda^t}_+$ the isomorphism that sends $(D_1,\ldots,D_n)$ to
$(D_n,\ldots,D_1)$.

\begin{Def} For any smooth $\Qlb$-sheaf $E$ on $X$ define the sheaf 
$E^{\lambda}_-$ on $X^{\lambda}_-$ by
$$
E^{\lambda}_-=\cH^{2a(\lambda)}(i_{\lambda}^*\cL^d_E)(a(\lambda))
$$
Define also the sheaf $E^{{\lambda}^t}_+$ on $X^{{\lambda}^t}_+$ by
$E^{{\lambda}^t}_+=t_*E^{\lambda}_-$. 
\end{Def}

 Let $S^{\lambda}\to X^{\lambda}_-$ be the vector bundle
whose fibre over $(D_1,\ldots,D_n)$ is the vector space of
collections $(\sigma_1,\ldots,\sigma_{n-1})$, where
$$
\sigma_i\in \Hom(\Omega^{n-i-1},\Omega^{n-i}(D_i)/\Omega^{n-i})
$$
By $\mu_S: S^{\lambda}\to \A^1$ we will denote the map
that at the level of $k$-points sends $(\sigma_1,\ldots,\sigma_{n-1})$
to the sum of $n-1$ classes in $k\iso \Ext^1(\Omega^{n-i-1},
\Omega^{n-i})$ corresponding to the pull-backs of
$$
0\to\Omega^{n-i}\to \Omega^{n-i}(D_i)\to\Omega^{n-i}(D_i)/
\Omega^{n-i}\to 0
$$
with respect to $\sigma_i: \Omega^{n-i-1}\to \Omega^{n-i}(D_i)/
\Omega^{n-i}$. 

Let $\cW^{\lambda}$ be the stack of collections:
$(D_1,\ldots,D_n)\in X^{\lambda}_-$ and a flag 
$(F^1\subset\ldots\subset F^n)$ of coherent torsion sheaves on $X$
with trivializations 
$$
F^i/F^{i-1}\iso
\Omega^{n-i}(D_i)/\Omega^{n-i}
$$
for $i=1,\ldots,n$. 
The projection $\tau:\cW^{\lambda}\to X^{\lambda}_-$ is a
generalized affine fibration of rank zero. 

 Let $\kappa:\cW^{\lambda}\to S^{\lambda}$ be the morphism over
$X^{\lambda}_-$ defined as follows. Given an $S$-point of $\cW^{\lambda}$,
consider for $i=1,\ldots,n-1$ the exact sequence
$$
0\to\Omega^{n-i-1}\to\Omega^{n-i-1}(D_{i+1})\to \Omega^{n-i-1}
(D_{i+1})/\Omega^{n-i-1}\to 0
$$
(Here $\Omega$ should be understood as the sheaf of relative differentials
$\Omega_{S\times X/S}$). It induces a map 
\begin{equation}
\label{morphism_1}
\HOM(\Omega^{n-i-1},\Omega^{n-i}(D_i)/\Omega^{n-i})
\to \EXT^1(\Omega^{n-i-1}(D_{i+1})/\Omega^{n-i-1},
\Omega^{n-i}(D_i)/\Omega^{n-i})
\end{equation}
which is an isomorphism of $\cO_{S\times X}$-modules, because $D_{i+1}\ge D_i$. 
The map $\kappa$ sends this point of $\cW^{\lambda}$ to 
$(\sigma_1,\ldots,\sigma_{n-1})$, where $\sigma_i$ is the global section of
$\HOM(\Omega^{n-i-1},\Omega^{n-i}(D_i)/\Omega^{n-i})$, whose image under
(\ref{morphism_1}) corresponds to the extension
$$
0\to F^i/F^{i-1}\to F^{i+1}/F^{i-1}\to F^{i+1}/F^i\to 0
$$

 Denote also by $\beta_{\cW}: \cW^{\lambda}\to {^{\le n}\Sh^d_0}$ the
morphism that sends $(F^1\subset\ldots\subset F^n)$ to $F^n$.

\begin{Pp}
\label{Pp_Cor_CS}
For any smooth $\Qlb$-sheaf $E$ on $X$ there is a
canonical isomorphism
\begin{equation}
\label{morphism_important}
\tau_!\, (\beta_{\cW}^*\, \cL^d_E\otimes
\kappa^*\mu_S^* \cL_{\psi})
\iso
E^{\lambda}_-
\end{equation}
\end{Pp}

 The proof is given in Sect. 3.3.

\bigskip
\noindent
3.2. \select{Local lemma}

\smallskip
\noindent
For the convenience of the reader, we begin with a local counterpart
of Proposition~\ref{Pp_Cor_CS}, working completely in a local
setting. 

  Let $\cO$ be a complete local $k$-algebra with residue field $k$,
which is regular of dimension one (that is, chosing a generator $\omega$
of the maximal ideal $\gm\subset\cO$, one identifies $\cO$ with
the ring $k[[\omega]]$ of formal power series of one
variable). Denote by $K$ the field of fractions of $\cO$. 
Let $\Omega$ be the completed module of relative differentials of
$\cO$ over $k$ (so, $\Omega$ is a free $\cO$-module generated by
$d\omega$). For $i\ge 0$ we write $\Omega^i$ for the $i$-th tensor power
of $\Omega$ (over $\cO$). 
For an integer $m$ denote by $\Omega^i(m)\subset
\Omega^i\otimes_{\cO} K$ the $\cO$-submodule generated by
$\omega^{-m}d\omega^{\otimes i}$. 

  Recall that we have fixed $\lambda\in
\Lambda_{n,d}^-$. Consider the stack $\cW^{\lambda}_{loc}$ classifying 
collections: a flag of torsion sheaves $(F^1\subset\ldots
\subset F^n)$ over $\Specf\cO$ with trivializations
$$
F^i/F^{i-1}\iso\Omega^{n-i}(\lambda_i)/\Omega^{n-i}
$$
for $i=1,\ldots,n$. (The subscript `loc' will stand for local counterparts
of certain stacks or morphisms).
Clearly, $\cW^{\lambda}_{loc}\to \Spec k$ is a generalized affine
fibration of rank zero. We also have the scheme $S^{\lambda}_{loc}$ whose set
of 
$k$-points is the
set of $(\sigma_1\ldots,\sigma_{n-1})$ with 
$$
\sigma_i\in\Hom(\Omega^{n-i-1},
\Omega^{n-i}(\lambda_i)/\Omega^{n-i})
$$
Besides, we have a map $(\mu_S)_{loc}: S^{\lambda}_{loc}\to\A^1$  
that at the level of $k$-points sends $(\sigma_1\ldots,\sigma_{n-1})$ to
$\sum \res \sigma_i$. One also defines a morphism
$\kappa_{loc}:\cW^{\lambda}_{loc}\to S^{\lambda}_{loc}$ in the same way as
$\kappa$. 

 Let $^{\le n}\Sh^d_0(\cO)$ be the stack classifying coherent torsion
sheaves $F$ on $\Specf\cO$ of length $d$ for which 
$\dim(F\otimes_{\cO}k)\le n$. 
It is stratified by locally closed substacks
$\Sh^{\nu}(\cO)$ indexed by $\nu\in \Lambda_{n,d}^+$. 
The stratum $\Sh^{\nu}(\cO)$ classifies sheaves
isomorphic to
$$
\cO/\gm^{\nu_1}\oplus\ldots\oplus\cO/\gm^{\nu_n}
$$
Let $\cB_{\nu}$ be the intersection cohomology sheaf associated to the
constant sheaf on the stratum $\Sh^{\nu}(\cO)$.
Let 
$$
\beta_{\cW,loc}: \cW^{\lambda}_{loc}\to {^{\le n}\Sh^d_0(\cO)}
$$
be the map sends $(F^1\subset\ldots\subset F^n)$ to $F^n$. 

 Local version of Proposition~\ref{Pp_Cor_CS} can be stated as follows.
\begin{Lm}
\label{Lm_1.11}
For any $\nu\in\Lambda_{n,d}^+$ we have canonically
$$
\R\Gamma_c( \cW^{\lambda}_{loc}, \;\;\beta_{\cW,loc}^*\,\cB_{\nu}\otimes
\kappa_{loc}^*(\mu_S)_{loc}^*\,\cL_{\psi})\iso
\left \{
\begin{array}{cc}
0, &  \mbox{ if } \nu^t\ne \lambda\\
\Qlb[-d](a(\lambda)),  & \mbox{ if }
\nu^t=\lambda
\end{array}
\right.
$$
\end{Lm}
\begin{Prf}
Denote by $\check{\cW}^{\lambda}_{loc}$ the stack of collections: a flag of
torsion sheaves $(\check{F}^1\subset\ldots\subset \check{F}^n)$ on
$\Specf\cO$ with trivializations
$$
\check{F}^i/\check{F}^{i-1}\iso\Omega^{n-i}/\Omega^{n-i}(-\lambda_{n-i+1})
$$
for $i=1,\ldots,n$.
We have an isomorphism $\cW^{\lambda}_{loc}\iso \check{\cW}^{\lambda}_{loc}$
that sends $(F^1\subset\ldots\subset F^n)$ to the flag
$(\check{F}^1\subset\ldots\subset \check{F}^n)$ with
$$
\check{F}^i=\EXT^1(F^n/F^{n-i}\, , \, \Omega^{n-1})
$$
for $i=1,\ldots,n$. This duality allows to switch between dominant and
anti-dominant weights of $\GL(n)$. 

 Put $\check{L}_i=\Omega^{n-1}\oplus\ldots\oplus\Omega^{n-i}$ for
$i=1,\ldots, n$. Denote by $\Gr^{d,+}(\check{L}_n)$ the moduli
scheme of $\cO$-sublattices $\check{\cR}\subset \check{L}_n$ such that 
$$
\dim(\check{L}_n/\check{\cR})=d
$$
Chosing a trivialization $\check{L}_n\iso\cO^n$, one identifies this scheme
with the connected component $\Gr^{d,+}=\Gr^{d,+}(\cO^n)$ of the positive
part of the affine grassmanian for $\GL(n)$. 

 We have a locally closed subscheme 
$\check{S}\hook{} \Gr^{d,+}(\check{L}_n)$ whose set of $k$-points 
consists of $\check{\cR}$
with the following property. If $\check{\cR}_i=\check{\cR}\cap \check{L}_i$
then the image of the inclusion
$\check{\cR}_i/\check{\cR}_{i-1}\hook{} \Omega^{n-i}$ is
$\Omega^{n-i}(-\lambda_{n-i+1})$ for $i=1,\ldots,n$.  

 We have a map $\,\check{\eta}_{loc}:\check{S}\to
\check{\cW}^{\lambda}_{loc}$  given by
$\check{F}^i=\check{L}_i/\check{\cR}_i$ for $i=1,\ldots,n$.   One checks
that $\check{\eta}_{loc}$ is an affine fibration of rank 
$a({\lambda}^t)$.
We also have a smooth and surjective map 
$$
\Gr^{d,+}(\check{L}_n)\to  {^{\le
n}\Sh^d_0(\cO)}
$$ 
that sends $\check{\cR}\subset \check{L}_n$ to
$\check{L}_n/\check{\cR}$, and we denote by
$\Gr^{\nu}(\check{L}_n)$ the preimage of the stratum
$\Sh^{\nu}(\cO)$ under this map. The Goresky-MacPherson extension of
$\Qlb[2\<\nu,\rho\>](\<\nu,\rho\>)$ from
$\Gr^{\nu}(\check{L}_n)$ to its closure is a perverse sheaf
denoted $\cA_{\nu}$ (cf. \cite{FGKV}). So, our assertion is nothing
else but the geometric Casselman-Shalika formulae (cf. \cite{N1,FGV,N2}):
$$
\R\Gamma_c(\check{S},\; \cA_{\nu}\otimes
\check{\eta}_{loc}^*\kappa_{loc}^* (\mu_S)_{loc}^*\cL_{\psi})\iso 
\left \{
\begin{array}{cc}
0, &  \mbox{ if } \nu^t\ne \lambda\\
\Qlb[-2\<\nu,\rho\>](-\<\nu,\rho\>),  & \mbox{ if }
\nu^t=\lambda
\end{array}
\right.
$$
\end{Prf}

\medskip

We will need the above lemma in a bit different form. 
Put $L_i=\Omega^{n-1}(\lambda_1)\oplus\ldots\oplus\Omega^{n-i}(\lambda_i)$
for $i=1,\ldots,n$. Let $\Gr^{d,+}(L_n)$ be the moduli scheme of
$\cO$-sublattices $\cR\subset L_n$ such that $\dim(L_n/\cR)=d$.
As in the proof of Lemma~\ref{Lm_1.11}, 
on $\Gr^{d,+}(L_n)$ we get a stratification by locally closed
subschemes $\Gr^{\nu}(L_n)$ indexed by $\nu\in\Lambda_{n,d}^+$,
and the perverse sheaves $\cA_{\nu}$.   

 By $S\subset \Gr^{d,+}(L_n)$ we denote the locally closed subscheme  whose
set of $k$-points consists of sublattices $\cR$ with the following property.
Let $\cR_i=\cR\cap L_i$. The condition is that the image of the inclusion
$$
\cR_i/\cR_{i-1}\hook{} L_i/L_{i-1}\iso \Omega^{n-i}(\lambda_i)
$$
is the sublattice $\Omega^{n-i}\subset \Omega^{n-i}(\lambda_i)$
for $i=1,\ldots,n$.  We also have a map $\eta_{loc}: S\to
\cW^{\lambda}_{loc}$ given by $F^i=L_i/\cR_i$ for $i=1,\ldots,n$. This is an
affine fibration of rank $a(\lambda)$. 

 Let $i:\Spec k \to S$ be the distinguished point that corresponds 
to $\cR=\Omega^{n-1}\oplus\ldots\oplus\cO$. Put $\mu_{loc}=
(\mu_S)_{loc}\comp\kappa_{loc}\comp\eta_{loc}$. 
We will use Lemma~\ref{Lm_1.11}
under the following form.

\begin{Lm} 
\label{Lm_1.12}
For any $\nu\in\Lambda_{n,d}^+$ we have canonically
$$
\R\Gamma_c(S,\;
\cA_{\nu}\otimes\mu_{loc}^*\cL_{\psi})\iso
\left \{
\begin{array}{cc}
0, &  \mbox{ if } \nu^t\ne \lambda\\
\Qlb[-2\<\lambda,\rho\>](-\<\lambda,\rho\>),  & \mbox{ if }
\nu^t=\lambda
\end{array}
\right.
$$
Moreover, this isomorphism is obtained by applying the
functor
$\R\Gamma_c$ to the composition of the canonical maps
$$
\cA_{\nu}\otimes\mu_{loc}^*\cL_{\psi}\to i_*
 i^*(\cA_{\nu}\otimes\mu_{loc}^*\cL_{\psi})
\to \tau_{\ge \; 2\<\lambda,\rho\>} ( i_*
i^*(\cA_{\nu}\otimes\mu_{loc}^*\cL_{\psi}))
$$
\end{Lm}
\begin{Prf}
As is easy to see, if $\nu^t\ne \lambda$ then the fibre 
$i^*(\cA_{\nu}\otimes\mu_{loc}^*\cL_{\psi})$
is placed in (usual) degrees strictly less then $2\<\lambda,\rho\>$.
For $\,\nu^t=\lambda$ this fibre equals 
$$
\Qlb[-2\<\lambda,\rho\>](-\<\lambda,\rho\>)
$$ 
Since the closure of $\Gr^{{\lambda}^t}(L_n)$ in $\Gr^{d,+}(L_n)$
is the union of strata
$\Gr^{\nu}(L_n)$ with $\nu\le {\lambda}^t$, 
our assertion follows from Lemma~\ref{Lm_1.11} combined with
the geometric statement due to B.C.Ngo (cf.~\cite{N2}, Lemma~5.2, p.14): 
\begin{itemize}
\item[i)] $S\cap \Gr^{\nu}(L_n)=\emptyset$ \ for \ $\nu< {\lambda}^t$
\item[ii)] $S\cap \Gr^{\nu}(L_n)$ is the point $\Spec k\hook{i} S$ \ for
\ $\nu={\lambda}^t$
\end{itemize}
\end{Prf}

\begin{Rem}
The shift in degree and the twist is calculated using the
following two formulas. For any weight $\nu$ of $\GL(n)$
we have $a(\nu)-a(\nu^t)=2\<\nu,\rho\>$ and
$2a(\nu)-2\<\nu,\rho\>=d(n-1)$,
where $d=\sum \nu_i$.
\end{Rem}

\bigskip
\noindent
3.3  \select{Proof of Proposition~\ref{Pp_Cor_CS}}

\smallskip\noindent
We return to our notation in the global case (as in Sect. 3.1). 
\Step 1
Denote by $\wt\cW^{\lambda}$ the scheme of collections: $(D_1,\ldots,D_n)
\in X^{\lambda}_-$, a digram
\begin{equation}
\label{diag_3}
\begin{array}{ccc}
L_1 & \subset\ldots\subset & L_n\\
\cup && \cup\\
\cR_1 & \subset\ldots\subset & \cR_n,
\end{array}
\end{equation}
where 
$
L_i=\Omega^{n-1}(D_1)\oplus\ldots\oplus\Omega^{n-i}(D_i)$ for
$i=1,\ldots,n$, and $(\cR_i)$ is a complete flag of vector subbundles on an
n-bundle $\cR_n$ such that the natural map
$$
\cR_i/\cR_{i-1}\hook{} L_i/L_{i-1}\iso\Omega^{n-i}(D_i)
$$
induces an isomorphism $\cR_i/\cR_{i-1}\iso \Omega^{n-i}$.
We have a map $\eta: \wt\cW^{\lambda}\to\cW^{\lambda}$ over
$X^{\lambda}_-$ given by 
$$
F^i=L_i/\cR_i
$$
for $i=1,\ldots,n$. This is an affine fibration of rank $a(\lambda)$,
so that $\eta_!\Qlb\iso \Qlb(-a(\lambda))[-2a(\lambda)]$. Put
$\tilde\tau=\tau\comp
\eta$. We will replace the functor $\tau_!(\cdot)$ by
$$
\tilde\tau_!\eta^*(\cdot)(a(\lambda))[2a(\lambda)]
$$
The advantage is that $\tilde\tau$ is representable whence $\tau$ is not.

 The morphism $\tilde\tau: \wt\cW^{\lambda}\to X^{\lambda}_-$ admits
a canonical section $\xi: X^{\lambda}_-\to \wt\cW^{\lambda}$ defined
by 
$$
\cR_i=\Omega^{n-1}\oplus\ldots\oplus\Omega^{n-i}
$$
for $i=1,\ldots,n$. Notice that $\xi$ is a closed
immersion. The following digram commutes
$$
\begin{array}{ccccc}
\cW^{\lambda} & \getsup{\eta} & \wt\cW^{\lambda} &
\toup{\tilde\tau} & X^{\lambda}_-\\
\downarrow\lefteqn{\scriptstyle\beta_{\cW}} &&
\uparrow\lefteqn{\scriptstyle{\xi}} &
\nearrow\lefteqn{\scriptstyle{\id}}\\
^{\le n}\Sh^d_0 & \getsup{i_{\lambda}} & X^{\lambda}_-
\end{array}
$$
Besides, the composition
$$
X^{\lambda}_-\toup{\xi}
\wt\cW^{\lambda}\toup{\eta}\cW^{\lambda}
\toup{\kappa} S^{\lambda}\toup{\mu_S}\A^1
$$
is the zero map. Now applying the functor $\tilde\tau_!$ to the
canonical morphism
$$
\eta^*\left(
\beta_{\cW}^*\cL^d_E\otimes
\kappa^*\mu_S^*\cL_{\psi} \right) \to
\xi_*\xi^*
\eta^*\left(
\beta_{\cW}^*\cL^d_E\otimes
\kappa^*\mu_S^*\cL_{\psi} \right) 
$$
we get a map 
\begin{equation}
\label{morphism_important_1}
\tau_!\left(
\beta_{\cW}^*\cL^d_E\otimes
\kappa^*\mu_S^*\cL_{\psi} \right) \to
i_{\lambda}^*\cL^d_E(a(\lambda))[2a(\lambda)]
\end{equation}
Define (\ref{morphism_important}) as the composition of
(\ref{morphism_important_1}) with the canonical map
$$
i_{\lambda}^*\cL^d_E(a(\lambda))[2a(\lambda)]
\to \cH^{2a(\lambda)}(i_{\lambda}^*\cL^d_E)(a(\lambda))
$$
Now we check that (\ref{morphism_important}) is an isomorphism fibre by
fibre.

\medskip
\Step 2
Fix a $k$-point $(D_1,\ldots,D_n)$ of
$X^{\lambda}_-$. Let $\,D_i=\sum_x \lambda_{i,x} \;x$. So,
the corresponding $\Lambda_n^-$-valued divisor on $X$ associates to
$x\in X$ the anti-dominant weight
$$
\lambda_x=(\lambda_{1,x},\ldots,\lambda_{n,x})\in
\Lambda_{n,d_x}^-  
$$
with $d_x=\sum_i \lambda_{i,x}$.  

For $i=1,\ldots,n$ put 
$L_i=\Omega^{n-1}(D_1)\oplus\ldots\oplus\Omega^{n-i}(D_i)$. 
For every closed point $x\in X$ let
$$
((L_1)_x\subset\ldots\subset(L_n)_x)
$$
be the restriction of the flag $(L_1\subset\ldots\subset
L_n)$ to $\Spec\hat\cO_{X,x}$. 

 Let $\Gr^{d_x,+}((L_n)_x)$ denote the moduli scheme of sublattices
$\cR\subset (L_n)_x$ such that 
$$
\dim((L_n)_x/\cR)=d_x
$$
By $S_x\subset \Gr^{d_x,+}((L_n)_x)$ we will
denote the locally closed subscheme whose set of $k$-points 
consists of sublattices $\cR\subset (L_n)_x$ with the
following property. 
Let $\cR_i=\cR\cap (L_i)_x$. The condition is that
the image of the natural inclusion
$$
\cR_i/\cR_{i-1}\hook{} (L_i)_x/(L_{i-1})_x\iso (\Omega^{n-i}(
\lambda_{i,x} \; x))_x
$$
is the sublattice $\Omega_x^{n-i}\subset (\Omega^{n-i}(
\lambda_{i,x} \; x))_x$ for $i=1,\ldots,n$. 
  
 For every $x\in X$ one defines the stack $\cW^{\lambda}_x$, the scheme 
$S^{\lambda}_x$ with morphisms 
$$
S_x\toup{\eta_x}\cW^{\lambda}_x\toup{\kappa_x} S^{\lambda}_x\toup{(\mu_S)_x}
\A^1,
$$
which are local counterparts of the corresponding stacks and morphisms 
for the weight $\lambda_x$. So, after the base change $(D_1,\ldots,D_n):
\Spec k\to X^{\lambda}_-$, the diagram
$$
\wt\cW^{\lambda}\toup{\eta}\cW^{\lambda}\toup{\kappa} S^{\lambda}
$$
becomes
$$
\prod_{x\in X}S_x\to\prod_{x\in X}\cW^{\lambda}_x\to\prod_{x\in X}
S^{\lambda}_x,
$$
the morphisms being the product of morphisms $\eta_x$ and $\kappa_x$
respectively. The restriction of $\mu_S: S^{\lambda}\to\A^1$ to $\prod_{x\in
X}S^{\lambda}_x$ is the sum of morphisms $(\mu_S)_x$.

 For $\nu\in\Lambda_{n,d_x}^+$
the perverse sheaf $\cA_{\nu}$ considered as a sheaf on 
$\Gr^{d_x,+}((L_n)_x)$ will be denoted 
$\cA_{\nu,x}\;$ (cf.  Sect. \hskip 0.4 em 3.2).  

  From (\cite{FGKV}, Proposition 3.1 and Lemma 4.2) it follows that the
restriction of $\eta^*\beta_{\cW}^*\cL^d_E$ to $\prod_{x\in X}
S_x$ is identified with $\boxtimes_{x\in X} F_x$, where $F_x$ is the
restriction of
$$
\oplusl_{\nu}
\cA_{\nu,x}[-d_x(n-1)](\frac{-d_x(n-1)}{2})\otimes E^{\nu}_x
$$
under the inclusion $S_x\hook{} \Gr^{d_x,+}((L_n)_x)$ (the sum being taken
over $\nu\in\Lambda_{n,d_x}^+$).

 Now combining Lemma~\ref{Lm_1.12} with (Theorem 3.3.8, 
\cite{L1}), we get the desired assertion. This
concludes the proof of Proposition~\ref{Pp_Cor_CS} $\;\square$.

\section{Geometric Whittaker models for $\GL(n)$}

4.1 \ \select{The stack $_n\cY_d$}

\medskip\noindent 
Consider the stack $_n\cX_d$ defined in Sect. 2.3. We impose
Pl\"ucker's  relations on a point $(L, (t_i))$ of $_n\cX_d$, which mean
that generically $(t_i)$ come from a complete flag of vector subbundles of
$L$. Our definition is justified by the following simple observation.

 Let $V$ be a vector space of dimension $n$ (over any field). 
For $n\ge k> i\ge 1$ let $\alpha_{k,i}:\wedge^k V\otimes\wedge^i V\to
\wedge^{k+1}V\otimes\wedge^{i-1}V$
be the contraction map that sends
$u\otimes(v_1\wedge v_2\wedge\ldots\wedge v_i)$ to 
$$
\sum_{j=1}^{i}
(-1)^j (u\wedge v_j)\otimes(v_1\wedge\ldots\wedge \hat 
{v}_j\wedge \ldots\wedge v_i)
$$

\begin{Lm}
Given nonzero elements $t_i\in \wedge^i V$ for $1\le i\le n$, the
following are equivalent:\\
1) There exists a complete flag of vector subspaces
$
0=V_0\subset V_1\subset\ldots\subset V_n=V
$
such that $t_i\in \wedge^i V_i\subset \wedge^i V$,\\
2) For $n\ge k> i\ge 1$ we have $\alpha_{k,i}(t_k\otimes
t_i)=0$. 
\end{Lm}
\begin{Prf}
The statement is obvious in characteristic zero. Let us give
an argument that holds in any characteristic. 
Write $e_1\wedge\ldots\wedge e_i$ for the image of 
$e_1\otimes\ldots\otimes e_i$ under $V^{\otimes i}\to \wedge^i V$.

 We construct by induction on $k$ the elements $e_1,\ldots,e_k\in V$ such that
$t_i=e_1\wedge\ldots\wedge e_i$ for $i=1,\ldots,k$. Let $e_1=t_1$, and
assume that $e_1,\ldots, e_{k-1}$ are already constructed.

 To construct $e_k$, we show by induction on $i$ that $t_k=e_1\wedge \ldots\wedge e_i\wedge
\omega_{k-i}$ for some $\omega_{k-i}\in\wedge^{k-i}V$, and define $e_k$ as $\omega_1$.

 First, since $\alpha_{k,1}(t_k\otimes t_1)=-t_k\wedge e_1=0$, we get $t_k=e_1\wedge 
\omega_{k-1}$ for some $\omega_{k-1}\in \wedge^{k-1}V$. Now assume that
$t_k=e_1\wedge\ldots\wedge e_{i-1}\wedge \omega_{k-i+1}$ for some $\omega_{k-i+1}\in
\wedge^{k-i+1}V$ with $i<k$. Then
$$
\alpha_{k,i}(t_k\otimes t_i)=\alpha_{k,i}(t_k\otimes (e_1\wedge\ldots\wedge e_i))= 
(-1)^i (t_k\wedge e_i)\otimes (e_1\wedge\ldots\wedge e_{i-1})=0
$$
It follows that 
$t_k\wedge e_i=0$. So, there exists $\omega_{k-i}\in\wedge^{k-i} V$ such that
$t_k=e_1\wedge \ldots\wedge e_i\wedge \omega_{k-i}$. We are done.
\end{Prf}

\medskip\smallskip

 Now  we define the closed substack $_n\cY_d\hook{}{_n\cX_d}$ by the
conditions $\alpha_{k,i}(t_k\otimes t_i)=0$ for $n\ge k> i\ge 1$, 
where
$$
\alpha_{k,i}:
\wedge^k L\otimes\wedge^i L\to \wedge^{k+1}L\otimes \wedge^{i-1}L
$$
are the contraction maps defined as above. Then the map $\varphi$ factors
through $_n\cQ_d\to {_n\cY_d}\hook{} {_n\cX_d}$. 

\smallskip

We stratify $_n\cY_d$ by locally closed substacks $\cV^{\lambda}_p\subset
{_n\cY_d}$ numbered by $\lambda\in\Lambda_{n,d}^p$. The stratum
$\cV^{\lambda}_p$ is defined by the condition: the degree of the divisor of zeros
of $t_i:\Omega^{(n-1)+\ldots+(n-i)}\hook{} \wedge^i L$ equals $\lambda_1+\ldots+\lambda_i$
for $i=1,\ldots,n$. Recall that
a point of $X^{\lambda}_p$ is a collection of divisors $(D_1,\ldots,D_n)$
on $X$ with $\deg(D_i)=\lambda_i$ and $D_1+\ldots+D_i\ge 0$ for all $i$.
So, the stack $\cV^{\lambda}_p$  classifies collections:
\begin{equation}
\label{collect_3}
(0=L'_0\subset L'_1\subset \ldots\subset L'_n=L, \; (s_i), \; (D_i)\in
X^{\lambda}_p),
\end{equation}
where $(L'_i)$ is a complete flag of subbundles on a rank $n$ vector
bundle $L$ and 
$$
s_i:\Omega^{(n-1)+\ldots+(n-i)}(D_1+\ldots+D_i)\iso \wedge^i
L'_i
$$
is an isomorphism. 

 Define the closed substack $\cV^{\lambda}\hook{}\cV^{\lambda}_p$ as
$\cV^{\lambda}_p\times_{X^{\lambda}_p} X^{\lambda}$. So, if $\lambda\notin
\Lambda_n$ then $\cV^{\lambda}$ is empty. 
Notice that the projection $\cV^{\lambda}_p\times_{_n\cY_d} \, {_n\cQ_d}\to
\cV^{\lambda}_p$ factors through $\cV^{\lambda}\hook{} \cV^{\lambda}_p$, and
for $\lambda\in\Lambda_{n,d}$ the corresponding morphism
$\cV^{\lambda}_p\times_{_n\cY_d} \, {_n\cQ_d}=
\cV^{\lambda}\times_{_n\cY_d} \, {_n\cQ_d}\to \cV^{\lambda}$ is an affine
fibration of rank $a(\lambda)$.

\bigskip
\noindent
4.2 \ \select{The sheaves $_n\cP^d_{E,\psi}$}

\begin{Def} For any smooth $\Qlb$-sheaf $E$ on $X$ put
$_n\cP^d_{E,\psi}=\varphi_!(_n\cF^d_{E,\psi})$. 
\end{Def}
 Clearly, the restriction of $_n\cP^d_{E,\psi}$ to a stratum
$\cV^{\lambda}_p$ of $_n\cY_d$ vanishes outside the closed substack
$\cV^{\lambda}$ of $\cV^{\lambda}_p$. For $\lambda\in\Lambda_{n,d}$
denote by $_n\cP^{\lambda}_{E,\psi}$ the restriction of $_n\cP^d_{E,\psi}$
to $\cV^{\lambda}$. 

\smallskip

 Define the closed substack $\cV^{\lambda}_-\hook{}\cV^{\lambda}$ as
$\cV^{\lambda}\times_{X^{\lambda}}X^{\lambda}_-$. Recall that the subscheme
$X^{\lambda}_-$ of $X^{\lambda}$ is given by the condition
$0\le\D_1\le\ldots\le D_n\,$, where $(D_i)\in X^{\lambda}$. Let 
$\mu_{\lambda}: \cV^{\lambda}_-\to \A^1$ be the map that at the
level of 
$k$-points sends (\ref{collect_3}) to the sum of $n-1$ classes in 
$$
k\iso \Ext^1(\Omega^{n-i-1}(D_i), \Omega^{n-i}(D_i))
$$
corresponding to the pull-backs of the successive extensions 
$$
\begin{array}{ccccc}
0\to & L'_i/L'_{i-1} & \to L'_{i+1}/L'_{i-1}\to  & L'_{i+1}/L'_i
& \to 0\\
     & \downarrow\lefteqn{\wr} && \downarrow\lefteqn{\wr}\\
     & \Omega^{n-i}(D_i) && \Omega^{n-i-1}(D_{i+1})
\end{array}
$$
with respect to the inclusion
$\Omega^{n-i-1}(D_i)\hook{}\Omega^{n-i-1}(D_{i+1})$. 

\begin{Pp}
\label{Pp_P_on_strata}
Let $E$ be a smooth $\Qlb$-sheaf on $X$ of rank $m$ and $\lambda\in
\Lambda_{n,d}$.  Then\\
1) $_n\cP^{\lambda}_{E,\psi}$ vanishes unless 
$$
\lambda_1=\ldots=\lambda_{n-m}=0 \eqno{(*)}
$$
2) Under the condition $(*)$ the complex
$_n\cP^{\lambda}_{E,\psi}$ is supported at $\cV^{\lambda}_-\hook{}
\cV^{\lambda}$, its restriction to $\cV^{\lambda}_-$ is 
isomorphic to the tensor product of
$$
\mu_{\lambda}^*\,\cL_{\psi}(\frac{b-2a(\lambda)}{2})[b-2a(\lambda)]
$$ 
with the inverse image
of $E^{\lambda}_-$ under $\cV^{\lambda}_-\to
X^{\lambda}_-$.
\end{Pp}

\begin{Rem}
\label{Rem_useful}
i) The sheaf $_n\cP^d_{E,\psi}$ was  also considered by Frenkel, Gaitsgory and Vilonen 
(\cite{FGV2}, 4.3). They show
that for any smooth $\Qlb$-sheaf $E$ on $X$, $_n\cP^d_{E,\psi}$ is a perverse sheaf
and the Goresky-MacPherson extension of its restriction to any nonempty open substack
of $_n\cY_d$. Besides, the Verdier dual of $_n\cP^d_{E,\psi}$ is
canonically isomorphic to $_n\cP^d_{E^*,\psi^{-1}}$ (4.6, 4.7, \select{loc.cit.}).
We notice that the stratification of $_n\cY_d$ used in (4.10, \select{loc.cit.})
is different from ours, so that our Proposition~\ref{Pp_P_on_strata} is a 
strengthened version of (4.13, \select{loc.cit.}). 

 According to (\cite{FGV,FGV2}), 
in the case $\rk E=n$ the perverse sheaf $_n\cP^d_{E,\psi}$ 
can be  thought of as a geometric counterpart of the
Whittaker function  canonically attached to $E$. 

\medskip
\noindent
ii) For $m>0$ let $_n^m\cY_d\subset {_n\cY_d}$ denote the open
substack given by the conditions: the image of $t_i$ is a \select{line
subbundle} in $\wedge^i L$ for $i=1,\ldots,n-m$. In particular,
$_n^m\cY_d={_n\cY_d}$ for $m\ge n$. 
Then 1) of Proposition~\ref{Pp_P_on_strata} claims  that
$_n\cP^d_{E,\psi}$ is the extension by zero of its restriction to
$_n^m\cY_d\subset {_n\cY_d}$.

\medskip
\noindent
iii) The relation between the sheaves $_n\cP^d_{E,\psi}$ for different $n$ is as follows.
Let $_n^m\cQ_d$ be the preimage of $_n^m\cY_d$ under
$\varphi:{_n\cQ_d}\to{_n\cY_d}$. So, $_n^m\cQ_d$ is the
open substack of $_n\cQ_d$ parametrizing collections (\ref{collect_1})
such that $L/L_{n-m}$ is locally free. 

 Denote by $_n^m\EXT$ the stack of collections
$(L_1\subset\ldots\subset L_{n-m+1}\subset L, (s_i))$, where
$L/L_{n-m}$ is a vector bundle on $X$ of rank $m$, and
$s_i:\Omega^{n-i}\iso L_i/L_{i-1}$ is an isomorphism
($i=1,\ldots,n-m+1$). Let $\mu_n^m: {_n^m\EXT}\to \A^1$ be the composition
${_n^m\EXT}\to {_{n-m+1}\cQ_0}\toup{\mu} \A^1$, where the first arrow is the map
that forgets $L$.  

 Let $_m\cM$ be the stack of pairs $(\Omega^{m-1}\hook{} L)$, where
$L$ is a vector bundle on $X$ of rank $m$. Taking the quotient by
$L_{n-m}$, we get a map $_n^m\EXT\to {_m\cM}$, which is a
generalized affine fibration. 

 For $1\le m\le n$ there is a commutative diagram
$$
\begin{array}{ccc}
_n^m\cQ_d & \iso & {_m\cQ_d}\times_{_m\cM} {_n^m\EXT}\\
\downarrow &&
\downarrow\lefteqn{\scriptstyle \varphi\times\id}\\
_n^m\cY_d & \iso  & {_m\cY_d}\times_{_m\cM}{_n^m\EXT};
\end{array}
$$
where the left vertical arrow is the restriction of $\varphi$.
So, the restriction of $_n\cP^d_{E,\psi}$ to $_n^m\cY_d$
is isomorphic to 
$$
_m\cP^d_{E,\psi}\boxtimes (\mu_n^m)^*\cL_{\psi}[b(m,d)-b(n,d)]
(\frac{b(m,d)-b(n,d)}{2})
$$
\end{Rem}

\medskip
\noindent
4.3 \ \select{The support of $_n\cP^{\lambda}_{E,\psi}$}

\medskip
\noindent
In this subsection we prove the next lemma.
\begin{Lm} 
\label{Lm_support_P_lambda}
The complex $_n\cP^{\lambda}_{E,\psi}$ vanishes outside
the closed substack $\cV^{\lambda}_-\hook{}\cV^{\lambda}$.
\end{Lm}

This may be derived from the geometric Casselman-Shalika formulae, but
we will give a direct proof. We start with the following sublemma.
Given $\,\lambda\in
\,\Lambda_n^-$ and $\nu\in\Lambda_n$ , denote by
$\;\cU_{\lambda}^{\nu}$ the stack of collections: $(D_1,\ldots,D_n)\in
X^{\lambda}_-$ ,
$(D'_1,\ldots,D'_n)\in X^{\nu}$, a diagram 
\begin{equation}
\label{diag_1}
\begin{array}{cccc}
L'_1& \subset \ldots\subset & L'_n\\
\cup && \cup\\
L_1 & \subset\ldots\subset & L_n\, ,
\end{array}
\end{equation}
where $(L_i)$ (resp., $(L'_i)$) is a complete flag of vector subbundles
on a rank $n$ vector bundle $L_n$ (resp., $L'_n$) on
$X$ with  trivializations 
$$
L'_i/L'_{i-1}\iso\Omega^{n-i}(D_i+D'_i)
$$ 
such that the image of the inclusion
$L_i/L_{i-1}\hook{}L'_i/L'_{i-1}\iso\Omega^{n-i}(D_i+D'_i)$
equals $\Omega^{n-i}(D_i)$ for $i=1,\ldots,n$. Let 
$$
\varphi_{\lambda}^{\nu}: \cU_{\lambda}^{\nu}\to 
(X^{\lambda}_-\times
X^{\nu})\times_{X^{\lambda+\nu}}\cV^{\lambda+\nu}
$$
be the map that forgets the flag $(L_i)$. Here $X^{\lambda}_-\times
X^{\nu}\to X^{\lambda+\nu}$ denotes the summation of
divisors. The map $\varphi_{\lambda}^{\nu}$ is an affine fibration
of rank $a(\nu)$.

\smallskip

 Let $X_{\lambda,\nu}\hook{}X^{\lambda}_-\times X^{\nu}$ be the
closed subscheme defined by
$$
D_i\ge D_{i-1}+D'_{i-1}
$$
for $i=2,\ldots,n$. The composition
$X_{\lambda,\nu}\hook{} X^{\lambda}_-\times X^{\nu}\to
X^{\lambda+\nu}$ factors through $X^{\lambda+\nu}_-
\hook{}X^{\lambda+\nu}$. 

\smallskip

 We also have a map
$\cU_{\lambda}^{\nu}\to \cV^{\lambda}_-$
that forgets $(L'_i)$ and $(D'_i)$. By abuse of notation, the
composition of this map with $\mu_{\lambda}:\cV^{\lambda}_-\to\A^1$ will
also be denoted $\mu_{\lambda}$.

\begin{Slm} 
\label{Slm_1}
The complex
$(\varphi^{\nu}_{\lambda})_! \mu_{\lambda}^*\cL_{\psi}$ is supported
at the closed substack $X_{\lambda,\nu}\times_{X^{\lambda+\nu}_-}
\cV^{\lambda+\nu}_-$ of $(X^{\lambda}_-\times
X^{\nu})\times_{X^{\lambda+\nu}}\cV^{\lambda+\nu}
$, and is isomorphic to the inverse image of
$$
\mu_{\lambda+\nu}^*\,\cL_{\psi}[-2a(\nu)](-a(\nu))
$$
from $\cV^{\lambda+\nu}_-$. 
\end{Slm}

\vskip -0.4 em

\begin{Prf} Let us decompose $\varphi^{\nu}_{\lambda}$ into two
affine fibrations
$\cU^{\nu}_{\lambda}\to \tilde\cU^{\nu}_{\lambda}
\toup{\tilde\varphi^{\nu}_{\lambda}} (X^{\lambda}_-\times
X^{\nu})\times_{X^{\lambda+\nu}}\cV^{\lambda+\nu}
$ defined as follows. Let $\tilde\cU^{\nu}_{\lambda}$ be the stack of
collections:
\begin{itemize}
\item $(D_i)\in X^{\lambda}_-\;$ , $\;(D'_i)\in X^{\nu}$
\item a complete flag of vector bundles $(L'_1\subset\ldots\subset L'_n)$ on
$X$ with trivializations
\begin{equation}
\label{trivialization_1}
L'_i/L'_{i-1}\iso\Omega^{n-i}(D_i+D'_i)
\end{equation}
for $i=1,\ldots,n$
\item for $i=1,\ldots,n-1$ diagrams
$$
\begin{array}{ccccccc}
0\to & L'_i/L'_{i-1} & \to & L'_{i+1}/L'_{i-1} & \to &  L'_{i+1}/L'_i &
\to 0 \\
& \cup&& \cup &&
\cup\\
0\to & \Omega^{n-i}(D_i) & \to & F_i & \to & \Omega^{n-i-1}(D_{i+1}) &
\to 0,
\end{array}
$$
where each row is an exact
sequence of $\cO_X$-modules, and both left and right vectical 
arrows are compatible with (\ref{trivialization_1}).
\end{itemize}
Now define the morphism $\cU^{\nu}_{\lambda}\to \tilde\cU^{\nu}_{\lambda}$ by
$F_i=L_{i+1}/L_{i-1}$ for $i=1,\ldots,n-1$, where $L_i$ are from diagram
(\ref{diag_1}). One checks that this is an
affine fibration of rank
$$
(n-2)\nu_1+(n-3)\nu_2+\ldots+\nu_{n-2}
$$
Define also $\tilde\varphi^{\nu}_{\lambda}$
as the map that forgets all $F_i$. This is an affine fibration of rank
$\nu_1+\ldots+\nu_{n-1}$. 

 Clearly, $\mu_{\lambda}:\cU^{\nu}_{\lambda}\to\A^1$ is constant along
the fibres of $\cU^{\nu}_{\lambda}\to\tilde\cU^{\nu}_{\lambda}$. So, it
suffices to prove the sublemma in the case $n=2$.

 In this case a fibre of $\varphi^{\nu}_{\lambda}$
is the affine space of maps $\xi: L_2/L_1\to L'_2/L_1$ such that the diagram
commutes
$$
\begin{array}{cccc}
0\to L'_1/L_1\to & L'_2/L_1 & \to & L'_2/L'_1\to 0\\
                 & \uparrow\lefteqn{\scriptstyle\xi} &
\nearrow\lefteqn{\scriptstyle i}\\
& L_2/L_1\, ,
\end{array}
$$
where $i$ is the canonical inclusion compatible with trivializations.
On this affine space we have a free and transitive action of
$\Hom(L_2/L_1,\, L'_1/L_1)$. The restriction of
$\mu_{\lambda}^*\cL_{\psi}$ to this affine space is a sheaf that changes
under the action of $\Hom(L_2/L_1, L'_1/L_1)$ by a local system, say
$\tilde\mu^*\cL_{\psi}$, where 
$$
\tilde\mu: \Hom(L_2/L_1, L'_1/L_1)\to k
$$
is the following linear functional. 
It associates to $s\in \Hom(L_2/L_1 ,\,
L'_1/L_1)$ the class of the pull-back  of
\begin{equation}
\label{sequence_1}
0\to L_1\to L'_1\to L'_1/L_1\to 0
\end{equation}
under the composition 
$\cO(D_1)\hook{}\cO(D_2)\iso L_2/L_1\toup{s}L'_1/L_1$.  The sequence
(\ref{sequence_1}) is just
$$
0\to \Omega(D_1)\to
\Omega(D_1+D'_1)\to \Omega(D_1+D'_1)/\Omega(D_1)\to 0
$$
So, $\tilde\mu=0$ if and only if $D_2\ge D_1+D'_1$. Besides, under this
condition the pull-back of
$$
0\to L'_1\to L'_2\to L'_2/L'_1\to 0
$$
under $\cO(D_1+D'_1)\hook{} \cO(D_2+D'_2)\iso L'_2/L'_1$ is identified 
(after tensoring by $\cO(-D'_1)$) with the pull back of 
$$
0\to L_1\to L_2\to L_2/L_1\to 0
$$ 
under $\cO(\D_1)\hook{}\cO(D_2)\iso L_2/L_1$. Our assertion follows.
\end{Prf}

\bigskip

For $m\ge 0$ and $\nu\in\Lambda_{m,d}$ denote by 
$\Fl^{\nu}$ the stack of flags
$(F^1\subset\ldots\subset F^m)$, where 
$F^i$ is a coherent torsion sheaf on $X$ with
$\deg(F^i/F^{i-1})=\nu_i\,$ for $i=1,\ldots, m$. Let 
$\div^{\nu}:\Fl^{\nu}\to X^{\nu}$ denote the composition
$$
\Fl^{\nu}\to\Sh^{\nu_1}_0\times\ldots\times\Sh^{\nu_m}_0
\,\toup{\div\times\ldots\times\div}\, X^{\nu}
$$
Set also
$_n\cQ^{\nu}={_n\cQ_d}\times_{\Sh^d_0} \Fl^{\nu}$, where
$\Fl^{\nu}\to\Sh^d_0$ sends $(F^1\subset \ldots\subset F^m)$
to $F^m$. 

 Denote by $_n^mJ_d$ the set of $n\times m$-matrices 
$e=(e^j_i)$ $(1\le i\le n,\; 1\le j\le m$) with $e^j_i\in\ZZ_+$,
$\sum_{i,j} e^j_i=d$. We have a map $h:
{_n^mJ_d}\to\Lambda_{n,d}\times \Lambda_{m,d}$ that sends 
$e$ to $(\lambda,\nu)$, where $\lambda_i=\sum_j e^j_i$ and
$\nu_j=\sum_i e^j_i$. For $e\in{_n^mJ_d}$ put $Y^e=\prod_{i,j}
X^{(e^j_i)}$. So, $Y^e$ classifies matrices of effective divisors
$(D^j_i)$ on $X$ such that $\deg(D^j_i)=e^j_i$. 

\smallskip

 For every $\lambda\in\Lambda_{n,d}$ the stack 
$\cV^{\lambda}\times_{_n\cY_d}
\, {_n\cQ^{\nu}}$ 
is stratified by locally closed substacks $\cQ^e\hook{}
\cV^{\lambda}\times_{_n\cY_d}\, {_n\cQ^{\nu}}$ 
indexed by $e\in {_n^m J_d}$ such that
$h(e)=(\lambda, \nu)$. The stratum $\cQ^e$ is the stack classifying
collections:
\begin{itemize}
\item a diagram

\vskip -2.5em
\begin{equation}
\label{diag_2}
\begin{array}{ccccccc}
\!\!\! L^m_1 & \subset & \!\!\! L^m_2 & \subset & \ldots &
\subset & \!\!\! L^m_n\\
\cup && \cup &&&& \cup\\
L^{m-1}_1 & \subset &  L^{m-1}_2 & \subset & \ldots &
\subset &  L^{m-1}_n\\
\cup && \cup &&&& \cup\\
\vdots && \vdots &&&& \vdots\\
\!\!\! L^0_1 & \subset & \!\!\! L^0_2 & \subset & \ldots &
\subset & \!\!\! L^0_n
\end{array}
\end{equation}
where $L^j_i$ is a vector bundle of rank $i$ on $X$, and
all the maps are inclusions of $\cO_X$-modules
\item a matrix $(D^j_i)\in Y^e$
\item isomorphisms 
$
L^m_i/ L^m_{i-1}\iso \Omega^{n-i}(D^1_i+\ldots+ D^m_i)
$
such that the image of the inclusion
$$
L^j_i/ L^j_{i-1}\hook{} L^m_i/ L^m_{i-1}\iso \Omega^{n-i}(D^1_i+\ldots+
D^m_i)
$$ 
equals $\Omega^{n-i}(D^1_i+\ldots+D^j_i)$ \ ($i=1,\ldots,n$; $\;
j=0,\ldots,m$)
\end{itemize}
We have a natural map
$$
\varphi^e: \cQ^e\to Y^e\times_{X^{\lambda}} \cV^{\lambda}
$$
that forgets all the rows in (\ref{diag_2}) except
the top one. (Here $Y^e\to X^{\lambda}$ sends $(D^j_i)$ to $(D_i)$
with $D_i=\sum_j D^j_i$). 
The morphism $\varphi^e$ is an affine fibration of rank $a(\lambda)$.

\smallskip

 Denote by $Y^e_-\hook{} Y^e$ the closed subscheme given by the
conditions:
\begin{itemize}
\item[$1'$)] for $i\le n-j$ we have $D^j_i=0$ 
\item[$2'$)] for $1\le j\le m-1\,$ and $\,2\le i\le n\,$ we have
$\, D_i^1+\ldots+D_i^j\ge D_{i-1}^1+\ldots+D_{i-1}^{j+1}$
\end{itemize}

 The composition $\cQ^e\hook{} {_n\cQ^{\nu}}\to {_n\cQ_d}\toup{\mu}\A^1$ 
will be denoted by $\mu_e$. 
\begin{Slm}
\label{Slm_difficult}
The copmlex 
$
(\varphi^e)_!\mu_e^*\cL_{\psi}
$
is supported at 
$Y^e_-\times_{X^{\lambda}_-}
\cV^{\lambda}_-\hook{} Y^e\times_{X^{\lambda}}\cV^{\lambda}$ and
is isomorphic to the inverse image of 
$$
\mu_{\lambda}^*\cL_{\psi}(-a(\lambda))[-2a(\lambda)]
$$
from $\cV^{\lambda}_-$. 
\end{Slm}
\begin{Prf} Apply Sublemma~\ref{Slm_1} $m$ times forgetting successively
the rows in diagram (\ref{diag_2}) starting from the
lowest one and moving up.
\end{Prf}

\bigskip
\begin{Prf}\select{of Lemma~\ref{Lm_support_P_lambda}}

\smallskip
\noindent
Since $\cL^d_E$ is a direct summand of the Springer sheaf
$\Spr^d_E$ (cf. Sect.~\ref{Sect_Springer}), it suffices to show that
the restriction of 
$$
\varphi_!(\beta^*\Spr^d_E\otimes\mu^*\cL_{\psi})
$$
to $\cV^{\lambda}$ vanishes outside $\cV^{\lambda}_-$. Put
$\nu=(1,\ldots,1)\in\Lambda_{d,d}$. The composition $_n\cQ^{\nu}\to{_n\cQ_d}
\toup{\mu}\A^1$ will also be denoted by $\mu$. 
By the projection formulae, we have to consider the direct image with
respect to the projection 
\begin{equation}
\label{morphism_2}
\cV^{\lambda}\times_{_n\cY_d}{_n\cQ^{\nu}}\to \cV^{\lambda}
\end{equation}
of $\pr_2^*\mu^*\cL_{\psi}$ tensored by some local system that comes
from $X^{\nu}$. The stack $\cV^{\lambda}\times_{_n\cY_d}{_n\cQ^{\nu}}$ 
is stratified by locally closed
substacks $\cQ^e$ indexed by $e\in {_n^dJ_d}$ such that $h(e)=(\lambda,
\nu)$. The restriction of (\ref{morphism_2}) to $\cQ^e$ can be
decomposed as
$$
\cQ^e\toup{\varphi^e} Y^e\times_{X^{\lambda}}\cV^{\lambda}\to
X^{\nu}\times_{X^{(d)}}\cV^{\lambda}\to \cV^{\lambda}
$$
So, our assertion follows from Sublemma~\ref{Slm_difficult}, because
the composition $Y^e_-\hook{} Y^e\to X^{\lambda}$ factors through
$X^{\lambda}_-\hook{} X^{\lambda}$.
\end{Prf}

\smallskip

\begin{Rem} Using Sublemma~\ref{Slm_difficult}, one may also check that
for any $\lambda\in\Lambda_{n,d}^-$ and any smooth $\Qlb$-sheaf $E$ on $X$,
$\, E^{\lambda}_-\, [\lambda_n]$ is a perverse sheaf on
$X^{\lambda}_-$. 
\end{Rem}

\medskip
\noindent
4.4 \ \select{Proof of Proposition~ \ref{Pp_P_on_strata}}

\smallskip
\noindent
Recall that $\cV^{\lambda}_-\times_{_n\cY_d}{_n\cQ_d}$
is the stack classifying collections: $(D_1,\ldots,D_n)\in
X^{\lambda}_-$ , a diagram
\begin{equation}
\label{diagram_17}
\begin{array}{ccccccc}
L'_1  & \subset \ldots \subset & L'_n\\
\cup && \cup\\
L_1 & \subset\ldots\subset & L_n\, ,
\end{array}
\end{equation}
where $(L'_i)$ (resp., $(L_i)$) is a complete flag of vector  subbundles
on a rank $n$ vector bundle $L'_n$ (resp., $L_n$) on $X$ with
trivializations
$$
L'_i/L'_{i-1}\iso \Omega^{n-i}(D_i)
$$
such that the image of the natural inclusion
$L_i/L_{i-1}\hook{}L'_i/L'_{i-1}\iso\Omega^{n-i}(D_i)$ equals $\Omega^{n-i}$
for $i=1,\ldots,n$. 
Denote by 
$$
\eta: \cV^{\lambda}_-\times_{_n\cY_d}{_n\cQ_d}\to
\cV^{\lambda}_-\times_{X^{\lambda}_-}
\cW^{\lambda}
$$ 
the morphism over $\cV^{\lambda}_-$ , whose composition with the
projection 
$\cV^{\lambda}_-\times_{X^{\lambda}_-}
\cW^{\lambda}\to \cW^{\lambda}$
sends (\ref{diagram_17}) to the flag
$
(L'_1/L_1\subset L'_2/L_2\subset\ldots\subset L'_n/L_n)
$. 
One checks that $\eta$ is a (representable) affine fibration of rank
$a(\lambda)$. Further, the composition
$$
\cV^{\lambda}_-\times_{_n\cY_d}{_n\cQ_d}\toup{\eta} 
\cV^{\lambda}_-\times_{X^{\lambda}_-}
\cW^{\lambda}\;\toup{\id\times\kappa}\;
\cV^{\lambda}_-\times_{X^{\lambda}_-} S^{\lambda}
\;\toup{\mu_{\lambda}\times \mu_S}\; \A^1\times\A^1\toup{{\rm sum}}\A^1
$$
coincides with the restriction of $\mu:{_n\cQ_d}\to \A^1$ to the
substack $\cV^{\lambda}_-\times_{_n\cY_d}{_n\cQ_d}\hook{} {_n\cQ_d}$. So, our
assertion follows from Proposition~\ref{Pp_Cor_CS}. $\square$

\section{Proof of Theorem~\ref{Th_A}}
\label{Sect_Prf_Th_A}

Recall the map $\phi: {_n\cX_d}\to X^{(d)}$ (cf. Sect.
2.3). By abuse of notation, its restriction to $_n\cY_d\hook{}
{_n\cX_d}$ is also denoted $\phi$. We let
$$
\tilde\pi: {_n\cY_d}\times_{_n\cM_d}{_n\cY_d}\to
X^{(d)}\times_{\Pic^d X} X^{(d)}
$$
be the morphism $\phi\times\phi$. By $\tilde\pi'$ we will denote the
restriction of
$\tilde\pi$ to the diagonal $_n\cY_d\hook{}
{_n\cY_d}\times_{_n\cM_d}{_n\cY_d}$. Clearly, Theorem~\ref{Th_A} is
equivalent to the fact that the natural map
$$
\tilde\pi_!({_n\cP^d_{E,\psi}}\boxtimes {_n\cP^d_{E',\psi^{-1}}})\to
\tilde\pi'_!({_n\cP^d_{E,\psi}}\otimes {_n\cP^d_{E',\psi^{-1}}})
$$
is an isomorphism. 
For $\lambda,\nu\in\Lambda_{n,d}$ we denote by 
$\tilde\pi^{\lambda,\nu}$ the restriction of $\tilde\pi$ to
the substack 
$$
\cV^{\lambda}\times_{_n\cM_d}\cV^{\nu}\hook{}
{_n\cY_d}\times_{_n\cM_d}{_n\cY_d}
$$ 
In the case $\lambda=\nu$ 
we write $(\tilde\pi^{\lambda,\lambda})'$ for the restriction of
$\tilde\pi^{\lambda,\lambda}$ to the diagonal $\cV^{\lambda}\hook{}
\cV^{\lambda}\times_{_n\cM_d}\cV^{\lambda}$.

 Using the stratification of
${_n\cY_d}\times_{_n\cM_d}{_n\cY_d}$ induced by both stratifications
of the first and the second multiple (cf. Sect.~4.1),
Theorem~\ref{Th_A} is reduced to the following statement.

\begin{Pp}
\label{Pp_on_Th_A}
For any $\lambda,\nu\in\Lambda_{n,d}$ 
the direct image
$(\tilde\pi^{\lambda,\nu})_!({_n\cP^{\lambda}_{E,\psi}}\boxtimes\,
{_n\cP^{\nu}_{E',\psi^{-1}}})$ 
vanishes unless $\lambda=\nu$. 
Under the condition $\lambda=\nu$ the natural map
$$
(\tilde\pi^{\lambda,\lambda})_!({_n\cP^{\lambda}_{E,\psi}}\boxtimes
{_n\cP^{\lambda}_{E', \psi^{-1}}})\to
(\tilde\pi^{\lambda,\lambda})'_!({_n\cP^{\lambda}_{E,\psi}}\otimes
{_n\cP^{\lambda}_{E',\psi^{-1}}})
$$
is an isomorphism.
\end{Pp}

\begin{Prf}
Put $\cV_1=\cV^{\lambda}_-\times_{_n\cM_d}\cV^{\nu}_-$. 
The restriction of $\tilde\pi^{\lambda,\nu}$ to $\cV_1$ can be
decomposed as
$$
\cV_1 \, \toup{^{1}\pi}\,
X^{\lambda}_-\times_{\Pic^d X} X^{\nu}_-
\to X^{(d)}\times_{\Pic^d X} X^{(d)},
$$ 
where $^{1}\pi$ is the product of two projections $\cV^{\lambda}_-\to
X^{\lambda}_-$ and $\cV^{\nu}_-\to X^{\nu}_-$.
In the case $\lambda=\nu$ we denote by
$\diag: \cV^{\lambda}_-\hook{} \cV_1$ the diagonal map.  

 By Proposition~\ref{Pp_P_on_strata}, our assertion
is reduced to the next lemma.

\begin{Lm}
For any $\lambda,\nu\in\Lambda_{n,d}^-$ 
the direct image
$\,{^{1}\pi_!}(\mu_{\lambda}^*\cL_{\psi}\boxtimes
\mu_{\nu}^*\cL_{\psi^{-1}})$ vanishes unless $\lambda=\nu$. Under the
condition $\lambda=\nu$ the natural map
$$
^{1}\pi_!(\mu_{\lambda}^*\cL_{\psi}\boxtimes
\mu_{\lambda}^*\cL_{\psi^{-1}})\to
{^{1}\pi_!}(\diag)_*(\diag)^*(\mu_{\lambda}^*\cL_{\psi}\boxtimes
\mu_{\lambda}^*\cL_{\psi^{-1}})
$$
is an isomorphism.
\end{Lm}

We need the next straightforward sublemma.
Given a divisor $D$ and a coherent sheaf $M$ on $X$ with a
section $\cO(D)\toup{s} M$, denote by
$\EXT_{M,D}$ the stack classifying extensions of
$\cO_X$-modules $0\to\Omega(D)\to ?\to M\to 0$, and by
$\mu_s:\EXT_{M,D}\to\A^1$ the map that sends this extension to the class of
its pull-back under $s$. 

\begin{Slm} 
\label{Slm_for_referee}
If $s\ne 0$ then $\RG_c(\EXT_{M,D},\;\mu_s^*\cL_{\psi})=0$.
$\square$
\end{Slm}

\begin{Prf}\select{of the lemma }
The stack $\cV_1$ classifies
collections: $(D_1,\ldots,D_n)\in X^{\lambda}_-$,
$(D'_1,\ldots,D'_n)\in X^{\nu}_-$, two flags
$(L_1\subset\ldots\subset L_n=L)$ and $(L'_1\subset\ldots\subset
L'_n=L)$ of subbundles on a rank $n$ vector bundle $L$ on $X$
with trivializations
$$
s_i: \Omega^{n-i}(D_i)\iso L_i/L_{i-1}\;\; \mbox{and}\;\;
s'_i: \Omega^{n-i}(D'_i)\iso L'_i/L'_{i-1}
$$ 
for $i=1,\ldots,n\,$
such that $(L_1, s_1)$ and $(L'_1,s'_1)$ coincide (in particular,
we have $D_1=D'_1$).

 Let $\cV_i\hook{} \cV_1$ be the closed
substack defined by the condition: the flags 
$$
(L_1\subset\ldots\subset
L_i\, , \, (s_j)_{j=1,\ldots,i})\;\; \mbox{and} \;\;
(L'_1\subset\ldots\subset L'_i\, , \, (s'_j)_{j=1,\ldots,i})
$$ 
coincide. Let $^i\cL_{\psi}$ be the restriction of 
$\mu_{\lambda}^*\cL_{\psi}\boxtimes \mu_{\nu}^*\cL_{\psi^{-1}}$ under
$\cV_i\hook{}\cV_1$. Let also 
$$
^{i}\pi: \cV_i\to X^{\lambda}_-\times_{\Pic^d X} X^{\nu}_-
$$
be the restriction of $^{1}\pi$ to $\cV_i$. Arguing by induction, we
will show that for every $i=1,\ldots,n-1$ the natural map
$$
(^{i}\pi)_!(^{i}\cL_{\psi})\to (^{i+1}\pi)_!(^{i+1}\cL_{\psi})
$$
is an isomorphism.

 To do so, denote by $\cN_i$ the stack of collections: 
$(D_1,\ldots,D_n)\in X^{\lambda}_-$,
$(D'_1,\ldots,D'_n)\in X^{\nu}_-$, two flags
$(M_{i+1}\subset\ldots\subset M_n=M)$ and
$(M'_{i+1}\subset\ldots\subset M'_n=M)$ of subbundles 
on a rank $n-i$ vector bundle $M$ on $X$ with trivializations
$$
s_j: \Omega^{n-j}(D_j)\iso M_j/M_{j-1}\;\; \mbox{and}\;\;
s'_j: \Omega^{n-j}(D'_j)\iso M'_j/M'_{j-1}
$$ 
for $j=i+1,\ldots,n$. 
Let also $'\cN_i\hook{}\cN_i$ be the closed substack defined by the
condition: 
$$
(M_{i+1},s_{i+1})\;\; \mbox{and}\;\; (M'_{i+1}, s'_{i+1})
$$
coincide. Taking the quotient by $L_i=L'_i$, we get a morphism
$\gamma:\cV_i \to \cN_i$, which is a generalized affine 
fibration. Further, we
have a commutative diagram
$$
\begin{array}{ccccc}
\cV_{i+1} & \hook{} & \cV_i \\
\downarrow && \downarrow\lefteqn{\scriptstyle{\gamma}} &
\searrow\lefteqn{\scriptstyle{ {^i\pi}}}\\ '\cN_i & \hook{} & \cN_i 
& \to &  X^{\lambda}_-\times_{\Pic^d X} X^{\nu}_-,
\end{array}
$$
where the square is cartesian. Applying Sublemma~\ref{Slm_for_referee}
for the section $s_{i+1}-s'_{i+1}:\Omega^{n-i-1}(D_i)\to M$, one 
checks that the complex $\gamma_!(^i\cL_{\psi})$ is supported at 
$'\cN_i$, and our assertion follows. 
\end{Prf}

\smallskip
\noindent
\end{Prf}(Proposition~\ref{Pp_on_Th_A})

\smallskip

 This concludes the proof of Theorem~\ref{Th_A}.

\section{Proof of Theorems~\ref{Th_B} and \ref{Th_C}}

6.1 \ \select{Plan of the proof}

\medskip\noindent
Proposition~\ref{Pp_P_on_strata} admits the following corolary.

\begin{Cor} 
\label{Cor_1}
For any smooth $\Qlb$-sheaves $E, E'$ on $X$ the complex
$$
_n\cS^d_{E,E'}\df
f_!(_n\cF^d_{E,\psi}\boxtimes{_n\cF^d_{E',\psi^{-1}}})(-d)[-2d]
$$
is a sheaf on $X^{(d)}$ placed in degree zero. It has a canonical filtration 
by constructable subsheaves such that
$$
\graded\;\, {_n\cS^d_{E, E'}}=\oplus_{\lambda\in\Lambda_{n,d}^+}
\; \pi^{\lambda}_* (E^{\lambda}_+\otimes
E^{\p\,\lambda}_{\, +})
$$
For each $r\le n$ there is a canonical inclusion $_r\cS^d_{E,E'}\subset {_n\cS^d_{E,E'}}$
compatible with filtrations.
\end{Cor}
\begin{Prf}
By the projection formulae, $\, _n\cS^d_{E,E'}\,\iso\,
\phi_!(_n\cP^d_{E,\psi}\otimes{_n\cP^d_{E',\psi^{-1}}})(-d)[-2d]$. Calculate
this direct image with respect to the stratification of $_n\cY_d$ by locally
closed substacks $\cV^{\lambda}_p$ indexed by $\lambda\in\Lambda_{n,d}^p$
(cf. Sect.~4.1). Since the natural map
$\cV^{\lambda}_-\to X^{\lambda}_-$ is a generalized affine fibration of rank
$b-d-2a(\lambda)$, our first assertion follows from
Proposition~\ref{Pp_P_on_strata}.

 Recall that we have open substacks $_n^r\cY_d\subset {_n\cY_d}$ for $r\le n$
(cf. iii) of Remark~\ref{Rem_useful}, Sect. 4.2). 
Let $^{\le r}\phi$ be the restriction of $\phi$ to $_n^r\cY_d$. By  
\select{loc.cit.},  
$^{\le r}\phi_!(_n\cP^d_{E,\psi}\otimes{_n\cP^d_{E',\psi^{-1}}})(-d)[-2d]\,\iso\, {_r\cS^d_{E,E'}}$. 
Our second assertion follows.
\end{Prf}

\medskip\smallskip

 This reduces our proof of Theorem~\ref{Th_B} to the
following steps. 
For $\nu\in\Lambda_{m,d}\, , \,\nu'\in\Lambda_{m',d}$ and $c\df(\nu,\nu')$
set $V^c=X^{\nu}\times_{X^{(d)}}X^{\nu'}$. Recall our notation
$_n\cQ^{\nu}= {_n\cQ_d}\times_{\Sh^d_0}\Fl^{\nu}$ (cf. Sect.~4.3). Let
$$
f^c: {_n\cQ^{\nu}}\times_{_n\cY_d}{_n\cQ^{\nu'}}\to V^c
$$
denote the composition ${_n\cQ^{\nu}}\times_{_n\cY_d}{_n\cQ^{\nu'}}\to
\Fl^{\nu}\times_{X^{(d)}}\Fl^{\nu'}\,\toup{\div^{\nu}\times\div^{\nu'}}\,
V^c$. The morphism $f^c$ is of finite type. 
Let also $^nf^c$ be the restriction of $f^c$ to the closed substack
$$
{_n\cQ^{\nu}}\times_{_n\cQ_d}{_n\cQ^{\nu'}}
\subset{_n\cQ^{\nu}}\times_{_n\cY_d}{_n\cQ^{\nu'}}
$$
First step is as follows. 
\begin{Pp}
\label{Pp_on_fc} The morphism $f^c$ is of relative dimension $\le b-d$, and
the natural map of the highest cohomology sheaves
$$
\R^{2(b-d)}(f^c)_!(\mu^*\cL_{\psi}\boxtimes\mu^*\cL_{\psi^{-1}})\to
\R^{2(b-d)}(^nf^c)_!\Qlb
$$
is an isomorphism.
\end{Pp}

 Further, set $W^c=\Fl^{\nu}\times_{\Sh^d_0}\Fl^{\nu'}$. Let $\div^c:
W^c\to V^c$ be the map $\div^{\nu}\times\div^{\nu'}$. Set 
$^{\le n}W^c=W^c\times_{\Sh^d_0}{^{\le n}\Sh^d_0}$. Let also
$^{\le n}\div^c$ denote the restriction of $\div^c$ to $\,{^{\le
n}W^c}\subset W^c$.  

 The morphism $^nf^c$ is decomposed as
${_n\cQ^{\nu}}\times_{_n\cQ_d}{_n\cQ^{\nu'}}\toup{\beta^c} {^{\le
n}W^c}\toup{^{\le n}\div^c} V^c$, where $\beta^c$ is the
natural projection. Second step is the next lemma.

\begin{Lm}
\label{Lm_relative_dimension}
i) $\beta^c$ is smooth and surjective with connected fibres of dimension
$b$.\\
ii) $\div^c$ is of relative dimension $\le -d$, so that 
$$
\R^{2(b-d)}(^nf^c)_!\Qlb(b-d)\iso\R^{-2d}(^{\le n}\div^c)_!\Qlb(-d)
$$
\end{Lm}

\smallskip

 Now assume $c=(\nu,\nu)$ with $\nu=(1,\ldots,1)\in\Lambda_{d,d}$. 
Let $_n\cS^c_{E,E'}$ denote the direct image under
$V^c\to X^{(d)}$ of the sheaf
$$
\R^{2(b-d)}(f^c)_!(\mu^*\cL_{\psi}\boxtimes\mu^*\cL_{\psi^{-1}})(b-d)
$$
tensored by the local system $(E^{\boxtimes\, d})\boxtimes
(E^{\p\,\boxtimes\, d})$ on $V^c$. The group $S_d\times S_d$ acts
naturally on $_n\cS^c_{E,E'}$. By Corolary~\ref{Cor_1}, 
$_n\cS^d_{E,E'}$ is the sheaf of $S_d\times S_d$-invariants of
$_n\cS^c_{E,E'}$. 

 Combining Lemma~\ref{Lm_relative_dimension} and Proposition~\ref{Pp_on_fc},
we learn that the complex $(^{\le n}\div)_!(\Spr^d_E\otimes\Spr^d_{E'})$ is
placed in degrees $\le -2d$, 
and there is a canonical $S_d\times S_d$-equivariant isomorphism
$$
_n\cS^c_{E,E'}\iso \R^{-2d}(^{\le
n}\div)_!(\Spr^d_E\otimes\Spr^d_{E'})(-d)
$$
Thus, Theorem~\ref{Th_B} is reduced to Theorem~\ref{Th_C}.

\bigskip
\noindent
6.2 \  \select{The stack $_i\tilde\cZ_d$}

\nopagebreak
\medskip
\noindent
For $0\le i\le n$ denote by $_i\tilde\cQ_d$ the stack classifying
collections
\begin{equation}
\label{collect_4}
(0=L_0\subset L_1\subset \ldots\subset L_i\subset F, \;\;\; (s_j)),
\end{equation}
where $F\in \Sh_i$, $(L_j)$ is a complete flag of vector subbundles
on a rank $i$ vector bundle $L_i$, $\deg(F/L_i)=d$, 
and $s_j:\Omega^{i-j}\iso L_j/L_{j-1}$ is an isomorphism 
$(j=1,\ldots, i)$. 

\vskip 1.5pt

 We have the open substack $_i\cQ_d\subset
{_i\tilde\cQ_d}$ 
given by the condition: $F$ is locally free. 
We also have a map $_i\tilde\cQ_d\to\Sh_i$ that sends (\ref{collect_4})
to $F$. Define a substack 
$$
_i\tilde\cZ_d\hook{}
{_i\tilde\cQ_d}\times_{\Sh_i}{_i\tilde\cQ_d}
$$
as follows. If $S$ is a scheme then an object 
\begin{equation}
\label{collect_5}
(F, (L_j, s_j), (L'_j, s'_j))
\end{equation}
of $\Hom(S, {_i\tilde\cQ_d}\times_{\Sh_i}{_i\tilde\cQ_d})$ 
lies in $\Hom(S, {_i\tilde\cZ_d})$ if the collections $(L_j,s_j)$ and
$(L'_j,s'_j)$ coincide outside a closed subscheme of $S\times X$
finite over $S$. 

\begin{Lm}
\label{Lm_closed_immersion}
 The map
$_i\tilde\cZ_d\hook{}
{_i\tilde\cQ_d}\times_{\Sh_i}{_i\tilde\cQ_d}$
is a closed immersion. In particular, the stack $_i\tilde\cZ_d$ is algebraic.
\end{Lm}
\begin{Prf}
An object (\ref{collect_5}) of $\Hom(S,
{_i\tilde\cQ_d}\times_{\Sh_i}{_i\tilde\cQ_d})$ gives rise
to a pair of sections
$$
t_j: \Omega^{(i-1)+\ldots+(i-j)}\iso
\wedge^j  L_j \to \wedge^j  F 
$$
and 
$$
t'_j: \Omega^{(i-1)+\ldots+(i-j)}\iso
\wedge^j  L'_j \to \wedge^j  F
$$
Clearly, (\ref{collect_5}) lies in $\Hom(S, {_i\tilde\cZ_d})$ if and only if
the support of $t_j-t'_j$ is a closed subscheme of $S\times X$ finite over
$S$ (for all $j=1,\ldots,i$). 

 Since $F$ is $S$-flat, $F$ (as well as its exterior powers) is locally free
outside some closed subscheme of $S\times X$ finite over $S$. So, our
assertion is a consequence of  the following sublemma, communicated to the
author by V.Drinfeld.

\begin{Slm}
\label{Lm_Drinfeld}
1) Let $F$ be any coherent sheaf on $S\times X$, which is locally free
outside a closed subscheme of $S\times X$ finite over $S$. Let $s$ be a
global section of $F$. Consider the following sub\-functor $Z$ of $S$ 
(on the category of $S$-schemes): a morphism $S'\to S$ belongs to
$Z(S')$ if the pull-back of $s$ to $S'\times X$ vanishes outside a closed
subscheme of $S'\times X$ finite over $S'$. Then the subfunctor $Z$ is 
closed. \\
2) Suppose in addition that $F$ is locally free. Then $S'\to S$ belongs to
$Z(S')$ if and only if the pull-back of $s$ to $S'\times X$ vanishes.
\end{Slm}
\begin{Prf}
If $r: S\times X\to S\times \P^1$ is a finite morphism over $S$ then the
functor $Z$ does not change if we replace $(X,F,s)$ by $(\P^1, r_*F, r_*s)$.
After localizing with respect to $S$ we may assume that $S$ is affine
and there is $r$ as above with $r_*F$ 
locally free over $S\times\A^1$. (Recall that if $S$ is noetherian then
any finite morphism $S\times X\to S\times \P^1$ over $S$ is flat (cf. SGA1, 
IV. 5.9)). So, we are reduced to the case $X=\P^1$ with
$F$ locally free over $S\times \A^1$. 

 If $S=\Spec R$ then we have the projective $R[t]$-module $M=\H^0(S\times
\A^1, F)$ and its element $s$. Represent $M$ as a direct summand of a free
$R[t]$-module $M'$. Clearly, $M'$ is also a free $R$-module. If $s_i\in R$
are the coordinates of $s\in M'$ (over $R$) then $Z$ is the closed subscheme
of $S$ defined by the equations $s_i=0$.  
\end{Prf}

\smallskip
\noindent
\end{Prf}(Lemma~\ref{Lm_closed_immersion})

\medskip

 We have an open substack $_i\cZ_d\subset {_i\tilde\cZ_d}$ given by the
condition: $F$ is locally free. In particular,
by 2) of Sublemma~\ref{Lm_Drinfeld}, $_n\cZ_d={_n\cQ_d}
\times_{_n\cY_d} \, {_n\cQ_d}$.
Besides, if (\ref{collect_5}) is a point of $_n\tilde\cZ_d$ then the sections
$$
\Omega^{(n-1)+\ldots+(n-n)}\iso\det L_n\hook{}\det F
$$
and 
$$
\Omega^{(n-1)+\ldots+(n-n)}\iso\det L'_n\hook{}\det F
$$
coincide, where $\det: \Sh_n\to\Pic X$ denotes 
the determinant map (cf. \cite{KM}).
 This yields a map $\tilde f: {_n\tilde\cZ_d}\to X^{(d)}$ whose
restriction to ${_n\cZ_d}$ coincides with $f$. We will see that
$\tilde f$ is of relative dimension $\le b-d$, but not
of finite type (the stack $_n\tilde\cZ_d$ even has infinitely many
irreducible components).  

For $k=0,\ldots,i$ we have a closed substack
$_i^k\tilde\cZ_d\hook{}{_i\tilde\cZ_d}$ given by the
condition: for a point (\ref{collect_5}) of $_i\tilde\cZ_d$ the flags
$$
(0=L_0\subset\ldots\subset L_k, (s_j)_{j=1,\ldots,k})
$$
and
$$
(0=L'_0\subset\ldots\subset L'_k, (s'_j)_{j=1,\ldots,k})
$$
coincide. Notice that taking the quotient by $L_k=L'_k$, 
one gets a map $_i^k\tilde\cZ_d\to
{_{i-k}\tilde\cZ_d}$, which is a generalized affine fibration. 
This observation will be a key point in the proof of
Proposition~\ref{Pp_on_fc}.

\bigskip
\noindent
6.3 \ \select{Dimensions counting}

\medskip
\begin{Prf}\select{of Lemma~\ref{Lm_relative_dimension}}

\smallskip\noindent
i) The map $\beta^c$ is obtained by base change from the map $\beta:
{_n\cQ_d}\to {^{\le n}\Sh^d_0}$, which is surjective and extends to a
generalized affine fibration ${_n\tilde\cQ_d}\to
\Sh^d_0$ that sends (\ref{collect_4}) to $F/L_n$.

\smallskip
\noindent
ii) We stratify $W^c$ by locally closed substacks
$\cU^e\hook{} W^c$ indexed by $e\in{_m^{m'}J_d}$ with $h(e)=(\nu,\nu')$. 
A point
\begin{equation}
\label{collect_6}
(F^1\subset\ldots\subset F^m=F, \;\; (F^1)'\subset\ldots\subset (F^{m'})'=F)
\end{equation}
of $W^c$ lies in $\cU^e$ if 
$$
\deg (F_i^j)= \sum_{k\le i, \;\; l\le j} e^l_k \mbox{\ \ \ \ for\ \ \ \ }
1\le i\le m, \;\; 1\le j\le m'
$$
where $F_i^j=F^i\cap (F^j)'$. If (\ref{collect_6}) is a point of $\cU^e$ then
for $1\le i\le m,\; 1\le j\le m'$ define 
$\tilde F_i^j\in \Sh_0$ from the
cocartesian square
$$
\begin{array}{ccc}
F^j_{i-1} & \to & \tilde F^j_i\\
\uparrow && \uparrow\\
F^{j-1}_{i-1} & \to  & F^{j-1}_i
\end{array}
$$
and put $G^j_i=F^j_i/\tilde F^j_i$. Set also
$$
\cW^e=\prod_{i,j}
\Sh^{e^j_i}_0
$$
The map $\cU^e\to \cW^e$ that sends (\ref{collect_6}) to the collection
$(G^j_i)$ is a generalized affine fibration of rank zero.
We have a map $\cW^e\to Y^e$ that sends $G^j_i$ to the collection
$(\div G^j_i)$, and define $\div^e: \cU^e\to Y^e$ as the composition
$\cU^e\to \cW^e\to Y^e$. Since for any $i\ge 0$ the morphism $\div: \Sh^i_0
\to X^{(i)}$ is of relative dimension $\le -i$, our assertion follows.
\end{Prf}

\medskip

 Define the stack $_i\tilde\cZ^c$ by the cartesian square
$$
\begin{array}{ccc}
_i\tilde\cZ^c & \to & {_i\tilde\cZ_d}\\
\downarrow && \downarrow\\
\Fl^{\nu}\times_{X^{(d)}}
\Fl^{\nu'} & \to & \Sh^d_0\times_{X^{(d)}}\Sh^d_0,
\end{array}
$$
where the right vertical arrow sends (\ref{collect_5}) to $(F/L_i\; ,
F/L'_i)$.
 Let $_i\cZ^c\subset {_i\tilde\cZ^c}$ denote the preimage of
$_i\cZ_d$ under $_i\tilde\cZ^c\to{_i\tilde\cZ_d}$.
In particular, we have 
$_n\cZ^c={_n\cQ^{\nu}}\times_{_n\cY_d}{_n\cQ^{\nu'}}$. Let
$$
\tilde f^c:{_n\tilde\cZ^c}\to V^c
$$
denote the composition ${_n\tilde\cZ^c}\to
\Fl^{\nu}\times_{X^{(d)}}\Fl^{\nu'}\toup{\div^{\nu}\times\div^{\nu'}} V^c$.
The restriction of $\tilde f^c$ to
${_n\cZ^c}$ coincides with $f^c$. Notice that $\tilde f^c$ is locally of
finite type, but not of finite type in general.

\begin{Lm} 
\label{Lm_b-d}
The map $\tilde f^c$ is of relative dimension $\le b-d$.
\end{Lm}
\begin{Prf}
\Step 1 The stack ${_n\cQ^{\nu}}\times_{_n\cY_d}{_n\cQ^{\nu'}}$ is stratified
by locally closed substacks $\cQ^e\times_{_n\cY_d}\cQ^{e'}$ indexed
by pairs $e\in{_n^m J_d}, e'\in{_n^{m'}J_d}$ such that there exists
$\lambda\in\Lambda_{n,d}$ with $h(e)=(\lambda,\nu)$, $h(e')=(\lambda,\nu')$.
(cf. Sect.~4.3).

 The restriction of $f^c:{_n\cZ^c}\to V^c$ to a stratum
$\cQ^e\times_{_n\cY_d}\cQ^{e'}$ is written as the composition
$$
\cQ^e\times_{_n\cY_d}\cQ^{e'}\,\toup{\;\varphi^e\times\;\varphi^{e'}}\, 
(Y^e\times_{X^{\lambda}}Y^{e'})\times_{X^{\lambda}}
\cV^{\lambda}\to Y^e\times_{X^{\lambda}}Y^{e'}
\to V^c
$$
Since $\varphi^e:\cQ^e\to Y^e\times_{X^{\lambda}}\cV^{\lambda}$ is an
affine fibration of rank $a(\lambda)$, and $\cV^{\lambda}\to X^{\lambda}$
is a generalized affine fibration of rank $b-d-2a(\lambda)$, it follows that
$f^c$ is of relative dimension $\le b-d$.

\smallskip
\Step 2 Stratify $\Sh_n$ by fixing the degree of the maximal torsion
subsheaf of $F\in\Sh_n$. Consider the induced stratification of
$_n\tilde\cZ^c$. A stratum $_n\tilde\cZ^c_k\subset {_n\tilde\cZ^c}$
classifies data: a point (\ref{collect_5}) of $_n\tilde\cZ_d$, two flags
of subsheaves 
\begin{equation}
\label{flag_Z_1}
(L_n\subset L_n^1\subset\ldots\subset L^m_n=F)
\end{equation}
and
\begin{equation}
\label{flag_Z_2}
(L'_n\subset (L_n^1)'\subset\ldots\subset (L_n^{m'})'=F),
\end{equation}
and an exact sequence 
$0\to F_0\to F\to M\to 0$ of $\cO_X$-modules, where
$F_0\in\Sh^k_0$ and $M$ is a vector bundle on $X$ of rank $i$.  

\smallskip

 The preimages of
flags (\ref{flag_Z_1}) and (\ref{flag_Z_2}) in $F_0$ give rise to a point
of $\Fl^{\lambda}\times_{\Sh^k_0} \Fl^{\lambda'}$ for some
$\lambda\in\Lambda_{m,k}$, $\lambda'\in\Lambda_{m',k}$. This yields
a stratification of $_n\tilde\cZ^c_k$ by locally closed substacks
$_n\tilde\cZ^c_{\lambda,\lambda'}\hook{}{_n\tilde\cZ^c_k}$ indexed by pairs
$\lambda\in\Lambda_{m,k}$, $\lambda'\in\Lambda_{m',k}$. 

\smallskip

 For an object of $_n\tilde\cZ^c_{\lambda,\lambda'}$ the vector bundle $M$
together with the images of the corresponding flags on $F$ defines a point
of $_n\cQ^{\nu-\lambda}\times_{_n\cY_{d-k}}{_n\cQ^{\nu'-\lambda'}}$. 
The natural forgetful map
$$
_n\tilde\cZ^c_{\lambda,\lambda'}\to (\Fl^{\lambda}\times_{\Sh^k_0}
\Fl^{\lambda'})\times 
({_n\cQ^{\nu-\lambda}}\times_{_n\cY_{d-k}}{_n\cQ^{\nu'-\lambda'}})
$$
is a generalized affine fibration of rank $nk$. Recall that $b=b(n,d)$
depends on $n$ and $d$ (cf. Sect.~2.2). By Step 1, 
$$
{_n\cQ^{\nu-\lambda}}\times_{_n\cY_{d-k}}{_n\cQ^{\nu'-\lambda'}}
\to X^{\nu-\lambda}\times_{X^{(d-k)}} X^{\nu'-\lambda'}
$$
is of relative dimension $\le b(n,d-k)-(d-k)=b(n,d)-d-nk+k$. By ii) of
Lemma~\ref{Lm_relative_dimension}, 
$$
\Fl^{\lambda}\times_{\Sh^k_0}
\Fl^{\lambda'}\to X^{\lambda}\times_{X^{(k)}} X^{\lambda'}
$$
is of relative dimension $\le -k$. Our assertion follows.
\end{Prf}

\bigskip
\noindent
6.4 \  \select{Proof of Proposition~\ref{Pp_on_fc}}

\medskip
\noindent
Let $_i^k\cZ_d$ be the preimage of $_i\cZ_d$  under
$_i^k\tilde\cZ_d\hook{} {_i\tilde\cZ_d}$. For $k=0,\ldots,i$ 
define the stacks $_i^k\cZ^c\subset {_i^k\tilde\cZ^c}$ 
by the cartesian squares
$$
\begin{array}{ccccc}
_i^k\cZ^c & \hook{} & {_i^k\tilde\cZ^c} & \hook{} & {_i\tilde\cZ^c}\\
\downarrow &&\downarrow && \downarrow\\
_i^k\cZ_d & \hook{} & {_i^k\tilde\cZ_d} & \hook{} & {_i\tilde\cZ_d}\\
\end{array}
$$
Denote by $^i\cL_{\psi}$ the restriction of 
$\mu^*\cL_{\psi}\boxtimes \mu^*\cL_{{\psi}^{-1}}$ under the composition
$$
^i_n\cZ^c\hook{}{_n\cZ^c}\to{_n\cZ_d}\,\iso\, {_n\cQ_d}
\times_{_n\cY_d} \, {_n\cQ_d}
$$ 
Let also $^if^c$ be the restriction of $f^c$ to
$^i_n\cZ^c\hook{}{_n\cZ^c}$.  
Arguing by induction, we will show that
the natural map
$$
R^{2(b-d)}(^if^c)_!(^i\cL_{\psi})\to
R^{2(b-d)}(^{i+1}f^c)_!(^{i+1}\cL_{\psi})
$$
is an isomorphism for $i=1,\ldots,n-1$. 

 The map $\mu:{_n\cQ_d}\to\A^1$ extends naturally to a morphism
$_n\tilde\cQ_d\to \A^1$ defined in the same way, it will also be denoted
by $\mu$. This allows to extend $^i\cL_{\psi}$ to
a local system $^i\tilde\cL_{\psi}$ on $_n^i\tilde\cZ^c$, where
$^i\tilde\cL_{\psi}$ is defined
as the retsriction of 
$\mu^*\cL_{\psi}\boxtimes
\mu^*\cL_{{\psi}^{-1}}$ under the composition
$$
_n^i\tilde\cZ^c\hook{}{_n\tilde\cZ^c}\to{_n\tilde\cZ_d}\hook{}
{_n\tilde\cQ_d}\times_{\Sh_i}{_n\tilde\cQ_d}
$$
We have the diagram
$$
\begin{array}{ccccc}
{}^{i+1}_n\cZ^c  & \hook{} & {{}^i_n\cZ^c} & \hook{\delta} &
{{}^i_n\tilde\cZ^c}\\
\downarrow & \searrow\lefteqn{\scriptstyle\alpha}&
\downarrow\lefteqn{\scriptstyle\beta} &
\swarrow\lefteqn{\scriptstyle\gamma}\\ 
{}^1_{n-i}\tilde\cZ^c & \hook{} &
{_{n-i}\tilde\cZ^c},
\end{array} 
$$
in which the square is cartesian, and $\gamma$ is a generalized
affine fibration of rank $b(n,d)-b(n-i,d)$. Since $\delta$ is
an open immersion, 
$$
\R^{top}\beta_!(^i\cL_{\psi})\to \R^{top}\gamma_!(^i\tilde\cL_{\psi})
$$
is an isomorphism over the image of (the smooth map) $\beta$. 
Applying Sublemma~\ref{Slm_for_referee} for $M=F/L_i$, $D=0$, and the section
$s_{i+1}-s'_{i+1}:\Omega^{n-i-1}\to M$, we learn that
$\R^{top}\gamma_!(^i\tilde\cL_{\psi})$ is supported at
${}^1_{n-i}\tilde\cZ^c$. Therefore,
$$
\R^{top}\beta_!(^i\cL_{\psi})\to \R^{top}\alpha_!(^{i+1}\cL_{\psi})
$$
is an isomorphism. So, $^if^c$ is decomposed as
$_n^i\cZ^c\toup{\beta} {_{n-i}\tilde\cZ^c}\to V^c$, where, by
Lemma~\ref{Lm_b-d}, the second map is of relative dimension $\le
b(n-i,d)-d$. Though $_{n-i}\tilde\cZ^c\to V^c$ is not of finite type in
general, its restriction to the image of $\beta$ is a morphism of finite
type. 

 This concludes the proof of Proposition~\ref{Pp_on_fc}.

\bigskip
\noindent
6.5 \  \select{Proof of Theorem~\ref{Th_C}}

\smallskip
\noindent
Recall that we have the map $h: {_m^{m'} J_d}\to
\Lambda_{m,d}\times \Lambda_{m',d}$, and for $e\in {_m^{m'}J_d}$ we write
$Y^e=\prod_{i,j} X^{(e_i^j)}$ (cf. Sect.~4.3). Let
$$
\norm: \mathop{\sqcup}\limits_{e\in h^{-1}(c)} Y^e \to V^c 
$$
be the map that sends a matrix $(D_i^j)\in Y^e$ to the collection $((D_i),
(D'_j))$, where $D_i=\sum_j D^j_i$ and $D'_j=\sum_i D^j_i$.

\begin{Lm}
i) The scheme $V^c$ is of pure dimension $d$, its
irreducible components are numbered by the set $h^{-1}(c)$. Namely, to
$e\in h^{-1}(c)$ there corresponds the component $\norm(Y^e)$.\\ 
ii) $\norm$ is the normalization of $V^c$ (more precisely, it is a finite
morphism, an isomorphism over an open dense subscheme of $V^c$, and
the scheme
$\;\mathop{\sqcup}\limits_{e\in h^{-1}(c)} Y^e$ is smooth). So,
$$
\norm_* \Qlb[d]\iso \IC,  
$$ 
where $\IC$ is the intersection cohomology sheaf on $V^c$.
\end{Lm}

\begin{Prf}
Stratify $V^c$ by locally closed subschemes $_e V^c\subset V^c$
indexed by $e\in h^{-1}(c)$. First, define $_e V^c$ as the open
subscheme of $Y^e$ given by the condition:
$$
\mbox{if } i>k, l>j \mbox{ then } D^j_i\cap D^l_k=\emptyset
$$
Then the composition $_e V^c\hook{} Y^e\toup{\norm} V^c$ is a locally
closed immersion. As a subscheme of $V^c$, $_e V^c$ is given by the 
condition: 
$$
\mbox{for all}\;\;\; i,j \;\;\;\mbox{we have} \;\;\;
\deg((\sum_{k\le i} D_k)\cap (\sum_{l\le j} D'_l))
=\sum_{k\le i,\;\;
l\le j} e^l_k
$$
For any $e\in h^{-1}(c)$ the scheme $_e V^c$ is smooth, nonempty and
irreducible of dimension $d$. This concludes the proof.
\end{Prf}

\begin{Lm}
\label{Lm_the_last_one}
There is a canonical isomorphism
$ R^{-2d}(\div^c)_! \Qlb (-d) \iso \norm_*\Qlb$.
\end{Lm}

 We need the following straightforward sublemma.

\begin{Slm}
\label{Slm_3}
i) Let $r:Y\to Y'$ be a separated morphism of schemes of finite
type. Assume that the fibres of $r$ are of dimension $\le d$. Let
$F$ be a smooth $\Qlb$-sheaf on $Y$, $U\subset Y$ be an open subscheme, 
and $r_U$ be the restriction of $r$ to $U$. Then the natural map
$\R^{2d}(r_U)_!F\to \R^{2d}r_!F$ is injective.

\smallskip
\noindent
ii) Let $(U^j)_{j\in J}$ be a stratification of $Y$ by locally
closed subschemes that comes from a filtration of $Y$ by closed
subschemes. Let $r^j$ be the restriction of $r$ to $U^j$. 
Then $\R^{2d}r_! F$ admits a filtration by
subsheaves with successive quotients being
$\R^{2d}r^j_!F$  $\;(j\in J)$. $\square$
\end{Slm}

\begin{Prf}\select{of Lemma~\ref{Lm_the_last_one}}

\smallskip\noindent
Recall that $W^c$ is stratified by locally closed substacks $\cU^e\hook{}W^c$
indexed by $e\in h^{-1}(c)$ and we have the maps $\div^e:\cU^e\to Y^e$
(cf. the proof of Lemma~\ref{Lm_relative_dimension}). The diagram commutes
$$
\begin{array}{ccc}
\cU^e & \hook{} & W^c\\
\downarrow\lefteqn{\scriptstyle{\div^e}} &&
\downarrow\lefteqn{\scriptstyle{\div^c}}\\
Y^e & \toup{\norm} & V^c
\end{array}
$$
We have $\R^{-2d}(\div^e)_!\Qlb(-d)\iso \Qlb$ 
canonically. Indeed, by K\"unneth formulae, this is reduced to the fact
that for any $i\ge 0$ the fibres of $\div:\Sh^i_0\to X^{(i)}$ are
connected of dimension $-i$. 

\smallskip

 By ii) of Sublemma~\ref{Slm_3}, on $R^{-2d}(\div^c)_! \Qlb (-d)$
there is a filtration parametrized by the set $h^{-1}(c)$ with successive
quotients being $(\norm^e)_*\Qlb$.  
We claim that any filtration with these successive quotients
degenerates canonically into a direct sum. Indeed, 
\begin{itemize}
\item[i)] the different successive quotients are supported on
different irreducible components of $V^c$, so our filtration
degenerates into a direct sum over some open dense subscheme of $V^c$;
\item[ii)] the sheaf $(\norm^e)_*\Qlb[d]$ is perverse, it is the
Goresky-MacPherson extension of its restriction to any open dense 
subscheme of $V^c$;
\item[iii)] the property ``perverse and the Goresky-MacPherson extension
of its restriction to a given open subscheme of $V^c$'' is preserved
for extensions. 
\end{itemize}
\end{Prf}

\pagebreak
\smallskip

 Finally, assume $c=(\nu,\nu)$ with $\nu=(1,\ldots,1)\in\Lambda_{d,d}$. 
Then the set $h^{-1}(c)$ is in natural bijection with $S_d$, and the map
$\norm$ becomes
$$
\mathop{\sqcup}\limits_{\sigma\in S_d} X^{\nu}_{\sigma} \to V^c,
$$
where $X^{\nu}_{\sigma}=X^{\nu}$, and $\norm$ sends a point
$(x_1,\ldots,x_d)\in X^{\nu}_{\sigma}$ to $((x_1,\ldots,x_d), (x_{\sigma 1},
\ldots, x_{\sigma d}))$. The action of $S_d\times S_d$ on $V^c$ 
lifts naturally to an action on
$$
(E^{\boxtimes\, d}\boxtimes
E^{\p\,\boxtimes\, d})\otimes \norm_*\Qlb,
$$ 
and it is easy to see that $\pr_!((E^{\boxtimes\,
d}\boxtimes E^{\p\,\boxtimes\, d})\otimes \norm_*\Qlb)^{S_d\times S_d}\iso
(E\otimes E')^{(d)}$ canonically, where $\pr: V^c\to X^{(d)}$ denotes the
projection.
 On the other hand, by Lemma~\ref{Lm_the_last_one}, 
$$
\R^{-2d}\div_!(\Spr^d_E\otimes\Spr^d_{E'})(-d)
\iso \pr_!((E^{\boxtimes\, d}\boxtimes
E^{\p\,\boxtimes\, d})\otimes \norm_*\Qlb)
$$
One checks that this isomorphism is $S_d\times S_d$-equivariant. Taking the 
invariants, one gets
$$
\R^{-2d}\div_!(\cL^d_E\otimes\cL^d_{E'})(-d)\iso
(E\otimes E')^{(d)}
$$

 By i) of Sublemma~\ref{Slm_3}, the natural map 
$\R^{-2d}(^{\le n}\div^c)_!\Qlb\to R^{-2d}(\div^c)_! \Qlb$ 
is an inclusion. It follows that
\begin{equation}
\label{fleche_1}
\R^{-2d}(^{\le n}\div)_!(\Spr^d_E\otimes\Spr^d_{E'})(-d)\to
\R^{-2d}\div_!(\Spr^d_E\otimes\Spr^d_{E'})(-d)
\end{equation}
is an inclusion. Taking the $S_d\times S_d$-invariants in (\ref{fleche_1}),
one gets an inclusion 
$_n\cS^d_{E,E'}\subset (E\otimes E')^{(d)}$, whose image is denoted
$^{\le n}(E\otimes E')^{(d)}$. Since $n$ and $d$ were arbitrary, 
Lemma~\ref{Lm_filtration} follows now from Corolary~\ref{Cor_1}, 
and the proof of Theorem~\ref{Th_C} is completed.

\smallskip

 So, Therem~\ref{Th_B} and Main Local Theorem are also proved.

\bigskip
\noindent
6.6  \select{Second proof of Theorem~\ref{Th_B}}

\medskip
\noindent
In this section we present an alternative proof of Theorem~\ref{Th_B} 
under the additional assumption: $\min\{\rk E,\rk E'\}\le n$. The idea
of this proof was suggested to the author by D.~Gaitsgory.

\smallskip 

Let $_n\Mod_d$ denote the stack classifying  modifications $(L\subset L')$ of 
rank $n$ vector bundles on $X$ with $\deg(L'/L)=d$. Let $q:{_n\Mod_d}\to\Sh^d_0$ be the map
that sends $(L\subset L')$ to $L'/L$, and $\supp: {_n\Mod_d}\to X^{(d)}$ denote
$\div\comp\, q$. For $d'\ge 0$
let $\gp_{\cY}:{_n\cY_d}\times_{\Bun_n}{_n\Mod_{d'}}\to {_n\cY_{d+d'}}$ be the map that
sends $((t_i), L\subset L')$ to $((t'_i), L')$, where $t'_i$ is the composition
$$
\Omega^{(n-1)+\ldots+(n-i)}\hook{t_i}\wedge^i L\hook{}\wedge^i L'
$$
The map $\gp_{\cY}$ is representable and proper. Let 
$\gq_{\cY}:{_n\cY_d}\times_{\Bun_n}{_n\Mod_{d'}}\to {_n\cY_d}$ denote the projection.
The map $\gq_{\cY}$ is smooth of relative dimension $nd'$. 

 The key ingredient is the Hecke property of Whittaker sheaves (\cite{FGV2}, 7.5). It
admits the following immediate corolary (the argument given in \select{loc.cit.} for
$\rk E=n$ holds, in fact, for $\rk E\le n$).
 
\begin{Pp} 
\label{Pp_Hecke_simplest}
For any  smooth $\Qlb$-sheaf $E$ on $X$  
and any $d\ge 0$ there is a  natural map
$$
(\gq_{\cY}\times\supp)_! \gp_{\cY}^*(_n\cP^{d+1}_{E,\psi})\to {_n\cP^d_{E,\psi}}\boxtimes
E (\frac{2-n}{2})[2-n], 
$$ 
which is an isomorphism if $\rk E\le n$. \QED
\end{Pp}

 Let $_n\wt{\Mod}_d$ be the stack of flags $(L_0\subset\ldots\subset L_d)$, where
$(L_i\subset L_{i+1})\in {_n\Mod_1}$ for all $i$. Let
$\wt{\supp}:{_n\wt{\Mod}_d}\to X^{d}$ be the map that sends $(L_0\subset\ldots\subset L_d)$ to
$(\div(L_1/L_0),\ldots,\div(L_d/L_{d-1}))$. Let
$p: {_n\cY_d}\times_{\Bun_n}{_n\wt{\Mod}_{d'}}\to {_n\cY_d}\times_{\Bun_n}{_n\Mod_{d'}}$ 
be the projection, and $\tilde\gp_{\cY}:
{_n\cY_d}\times_{\Bun_n}{_n\wt{\Mod}_{d'}}\to {_n\cY_{d+d'}}$ be the composition
$\gp_{\cY}\comp\, p$.

\begin{Cor}
\label{Cor_Hecke}
For any smooth $\Qlb$-sheaf $E$ on $X$ and any $d,d'\ge 0$ there is a natural map
\begin{equation}
\label{map_Cor_Hecke}
(\gq_{\cY}\times\wt{\supp})_!\, \tilde\gp_{\cY}^*(_n\cP^{d+d'}_{E,\psi})\to 
{_n\cP^d_{E,\psi}}\boxtimes E^{\boxtimes d'}(\frac{2d'-nd'}{2})[2d'-nd'],
\end{equation}
which is an isomorphism if $\rk E\le n$. 
\end{Cor}
\begin{Prf}
The map (\ref{map_Cor_Hecke}) is defined as follows. Let $\gp_{\cQ}: {_n\cQ_d}\times_{\Bun_n}
{_n\Mod_{d'}}\to {_n\cQ_{d+d'}}$ be the map that sends $(L_1\subset\ldots\subset L_n\subset 
L\subset L')$ to $(L_1\subset\ldots\subset L_n\subset L')$. Let $\tilde\gp_{\cQ}:
{_n\cQ_d}\times_{\Bun_n}{_n\wt{\Mod}_{d'}}\to {_n\cQ_{d+d'}}$ denote the composition
$$
{_n\cQ_d}\times_{\Bun_n}{_n\wt{\Mod}_{d'}}\to {_n\cQ_d}\times_{\Bun_n}{_n\Mod_{d'}}
\toup{\gp_{\cQ}} {_n\cQ_{d+d'}}, 
$$
where the first arrow is the projection. Consider the commutative diagram
$$
\begin{array}{ccc}
_n\cQ_d\times_{\Bun_n}{_n\wt{\Mod}_{d'}} & \toup{\tilde\gp_{\cQ}} & {_n\cQ_{d+d'}}\\
\downarrow\lefteqn{\scriptstyle \varphi\times\id} && \downarrow\lefteqn{\scriptstyle \varphi}\\
{_n\cY_d}\times_{\Bun_n}{_n\wt{\Mod}_{d'}} & \toup{\tilde\gp_{\cY}} & 
{_n\cY_{d+d'}}
\end{array}
$$
Since $_n\cF^{d+d'}_{E,\psi}$ is a direct summand of
$$
(\tilde\gp_{\cQ})_!(_n\cF^d_{E,\psi}\boxtimes \wt{\supp}^*E^{\boxtimes\, d'})[nd']
(\frac{nd'}{2}),
$$
it follows that $_n\cP^{d+d'}_{E,\psi}$ is a direct summand of 
$
(\tilde\gp_{\cY})_!(_n\cP^d_{E,\psi}\boxtimes \wt{\supp}^*E^{\boxtimes\, d'})[nd']
(\frac{nd'}{2}) 
$.
This yields a morphism 
$$
\tilde\gp_{\cY}^* (_n\cP^{d+d'}_{E,\psi})\to {_n\cP^d_{E,\psi}}\boxtimes \wt{\supp}^*
E^{\boxtimes \, d'}[nd'](\frac{nd'}{2})
$$
Since $\gq_{\cY}\times\wt{\supp}$ is smooth of relative dimension $d'(n-1)$, the desired map
is obtained from the last one by the adjointness.
To show that (\ref{map_Cor_Hecke}) is an isomorphism under the condition $\rk E\le n$, 
apply $d'$ times Proposition~\ref{Pp_Hecke_simplest}.
\end{Prf}

\medskip

Denote by $^{rss}X^{d'}\subset X^{d'}$ the open subscheme that parametrizes pairwise different
points $(x_1,\ldots,x_{d'})\in X^{d'}$ (here `rss' stands for `regular semisimple'). Let
$^{rss}_n\wt{\Mod}_{d'}$ be the preimage of $^{rss}X^{d'}$ under $\wt{\supp}$. The symmetric group
$S_{d'}$ acts on $^{rss}_n\wt{\Mod}_{d'}$, and the restriction of $\tilde\gp_{\cY}$
to ${_n\cY_d}\times_{\Bun_n}{_n^{rss}\wt{\Mod}_{d'}}$ is $S_{d'}$-invariant. So, the
action of $S_{d'}$ on ${_n\cY_d}\times_{\Bun_n}{_n^{rss}\wt{\Mod}_{d'}}$ 
lifts to an action on 
$$
\tilde\gp_{\cY}^*(_n\cP^{d+d'}_{E,\psi})
$$
Since the restriction $^{rss}_n\wt{\Mod}_{d'}\to {^{rss}X^{d'}}$ of $\wt{\supp}$ is 
$S_{d'}$-equivariant, $S_{d'}$ acts on the complex
$$
(\gq_{\cY}\times\wt{\supp})_!\tilde\gp_{\cY}^*(_n\cP^{d+d'}_{E,\psi})
$$
restricted to $_n\cY_d\times {^{rss}X^{d'}}$. On the other hand, $S_{d'}$ acts on $E^{\boxtimes 
\, d'}$ and, hence, on the right hand side of (\ref{map_Cor_Hecke}). Using the explicit 
description of the map (\ref{map_Cor_Hecke}) one easily proves the next lemma.

\begin{Lm}
\label{Lm_equiva}
 For any smooth $\Qlb$-sheaf $E$ on $X$ the map (\ref{map_Cor_Hecke})
restricted to $_n\cY_d\times {^{rss}X^{d'}}$ is $S_{d'}$-equivariant. \QED
\end{Lm}

\smallskip

 Recall that for any smooth $\Qlb$-sheaf $E$ on $X$ the Verdier dual of $_n\cP^d_{E,\psi}$
is canonically isomorphic to $_n\cP^d_{E^*,\psi^{-1}}$ \  (\cite{FGV2}, 4.7).
So, to prove Theorem~\ref{Th_B} we must establish a canonical isomorphism 
$$
\phi_*\R\HOM(_n\cP^d_{E,\psi}\, ,\; {_n\cP^d_{E',\psi}}) \iso
\HOM(E,E')^{(d)},
$$
where $\phi: {_n\cY_d}\to X^{(d)}$ is the map defined in Sect. 2.3. The statement 
of Theorem~\ref{Th_B} being symmetric with respect to interchanging $E$ and 
$E'$, we assume $\rk E\le n$. 

 For any smooth $\Qlb$-sheaf $E$ on $X$ set $_n\tilde{\cP}^d_{E,\psi}=\varphi_!(\beta^*\Spr^d_E
\otimes\mu^*\cL_{\psi})[b](\frac{b}{2})$. In other words, $_n\tilde{\cP}^d_{E,\psi}$ is
a complex on $_n\cY_d$ obtained by replacing in the definition of $_n\cP^d_{E,\psi}$ 
 Laumon's sheaf $\cL^d_E$ by Springer's sheaf $\Spr^d_E$. Theorem~\ref{Th_B} follows 
now from the next statement.

\begin{Pp}
\label{Pp_really_last}
Let $E,E'$ be smooth $\Qlb$-sheaves on $X$ with $\rk E\le n$. Then there exists a canonical 
$S_d$-equivariant isomorphism
$$
\phi_*\R\HOM(_n\cP^d_{E,\psi} ,\; {_n\tilde{\cP}^d_{E',\psi}}) \iso
\sym_*(\HOM(E,E')^{\boxtimes\, d})
$$
\end{Pp}
\begin{Prf} The idea is that Proposition~\ref{Pp_really_last} is a tautological consequence
of Corolary~\ref{Cor_Hecke} obtained by applying
the formalism of six functors. The equivariance property follows from Lemma~\ref{Lm_equiva}. 
The precise argument is as follows.

Consider the commutative diagram
$$
\begin{array}{ccccc}
 && {_n\cY_0}\times_{\Bun_n}{_n\Mod_d} & \toup{\gp_{\cY}} & _n\cY_d\\ 
& \swarrow\lefteqn{\scriptstyle \gq_{\cY}} & \downarrow\lefteqn{\scriptstyle 
\gq_{\cY}\times\supp} && \downarrow\lefteqn{\scriptstyle \phi}\\
_n\cY_0 & \gets & _n\cY_0\times X^{(d)} & \to & X^{(d)}
\end{array}
$$
Set for brevity $\Psi^0={_n\cP^0_{E,\psi}}$ (it does not depend on $E$, though does depend on 
$\psi$). By definition, 
$$
_n\tilde{\cP}^d_{E',\psi}\;\iso\; \gp_{\cY\, *}(\gq_{\cY}^*\Psi^0\otimes 
q^*\Spr^d_{E'})[nd](\frac{nd}{2})
$$
\begin{Lm}
$\;\gq_{\cY}^*\Psi^0\otimes q^*\Spr^d_{E'}\;\iso\; \R\HOM(q^*\Spr^d_{{E'}^*},\;
\gq_{\cY}^*\Psi^0)$
canonically and $S_d$-equivariantly. 
\end{Lm}
\noindent\select{Proof}
\ Using the fact that both $\gq_{\cY}$ and $\gq_{\cY}\comp p$ are smooth of relative dimension
$nd$, we get
\begin{multline*}
\gq_{\cY}^*\Psi^0\otimes q^*\Spr^d_{E'}\;\iso\; p_!(\wt{supp}^*{E'}^{\,\boxtimes\, d}\otimes
p^*\gq_{\cY}^*\Psi^0)\iso \;
p_*\R\HOM(\wt{\supp}^*({E'}^*)^{\boxtimes\, d},\; 
p^*\gq_{\cY}^*\Psi^0)\iso\\ 
p_*\R\HOM(\wt{\supp}^*({E'}^*)^{\boxtimes\, d},\; 
p^!\gq_{\cY}^!\Psi^0[-2nd](-nd))\;\iso\;
\R\HOM(p_!\wt{\supp}^*({E'}^*)^{\boxtimes\, d},\; \gq_{\cY}^*\Psi^0) \;\;\;\;\;\; \square
\end{multline*}

\medskip
\noindent
Using the above lemma, we get   
\begin{multline*}
\R\HOM(_n\cP^d_{E,\psi} ,\; {_n\tilde{\cP}^d_{E',\psi}})\; \iso\;
\gp_{\cY\, *}\R\HOM(\gp_{\cY}^*(_n\cP^d_{E,\psi}),\;  \gq_{\cY}^*\Psi^0\otimes 
q^*\Spr^d_{E'})[nd](\frac{nd}{2})\;\iso\\
\iso\; \gp_{\cY\, *}\R\HOM(\gp_{\cY}^*(_n\cP^d_{E,\psi})\otimes q^*\Spr^d_{{E'}^*},\; 
\gq_{\cY}^!\Psi^0)[-nd](\frac{-nd}{2})
\end{multline*}  

Let $j:{_n\cQ_0}\hook{}{_n\cY_0}$ denote the natural open immersion. Since $_n\cQ_0\to\Spec k$
is a generalized affine fibration, we have $\RG(_n\cQ_0, \;\Qlb)\iso\Qlb$. So, our assertion
is reduced to the next lemma.

\begin{Lm}
\label{Lm_13}
 There is a canonical $S_d$-equivariant isomorphism over $_n\cY_0\times X^{(d)}$
$$
(\gq_{\cY}\times\supp)_*\R\HOM(\gp_{\cY}^*(_n\cP^d_{E,\psi})\otimes q^*\Spr^d_{{E'}^*},
\; \gq_{\cY}^!\Psi^0)\;\iso\;
(j\times\id)_*(\Qlb\boxtimes \sym_*\HOM(E,E')^{\boxtimes\, d})[nd](\frac{nd}{2})
$$
\end{Lm}
\begin{Prf}
Let $\pr_{\cY}: {_n\cY_0}\times X^d\to {_n\cY_0}$ denote the projection.
Using the commutative diagram
$$
\begin{array}{ccc}
{_n\cY_0}\times_{\Bun_n}{_n\Mod_d} & \getsup{p} &  {_n\cY_0}\times_{\Bun_n}{_n\wt{\Mod}_d}\\
\downarrow\lefteqn{\scriptstyle \gq_{\cY}\times\supp} && \downarrow\lefteqn{\scriptstyle 
\gq_{\cY}\times\wt{\supp}}\\
_n\cY_0\times X^{(d)} & \getsup{\id\times\sym} & _n\cY_0\times X^d,
\end{array}
$$
we obtain
\begin{multline*}
(\gq_{\cY}\times\supp)_*\R\HOM(\gp_{\cY}^*(_n\cP^d_{E,\psi})\otimes q^*\Spr^d_{{E'}^*},
\; \gq_{\cY}^!\Psi^0)\;\iso\\
\hskip -12.6em
(\gq_{\cY}\times\supp)_*\R\HOM(p_!(\tilde{\gp}_{\cY}^*(_n\cP^d_{E,\psi})\otimes
\wt{\supp}^*({E'}^*)^{\boxtimes \, d}), \; \gq_{\cY}^!\Psi^0)\;\iso\\
\hskip -12.5em
(\gq_{\cY}\times\supp)_*p_*\R\HOM(\tilde{\gp}_{\cY}^*(_n\cP^d_{E,\psi})\otimes
\wt{\supp}^*({E'}^*)^{\boxtimes \, d}, \; p^!\gq_{\cY}^!\Psi^0)\;\iso\\
\hskip -3.2em
(\id\times\sym)_*(\gq_{\cY}\times\wt{\supp})_*\R\HOM(\tilde{\gp}_{\cY}^*(_n\cP^d_{E,\psi})\otimes
\wt{\supp}^*({E'}^*)^{\boxtimes \, d}, (\gq_{\cY}\times\wt{\supp})^!
\pr_{\cY}^!\Psi^0)\;\iso\\
\hskip -7.8em
(\id\times\sym)_*\R\HOM((\gq_{\cY}\times\wt{\supp})_!(\tilde{\gp}_{\cY}^*(_n\cP^d_{E,\psi})\otimes
\wt{\supp}^*({E'}^*)^{\boxtimes \, d}), \; \pr_{\cY}^!\Psi^0)\;\iso\\
\hskip -10.1em
(\id\times\sym)_*\R\HOM(\Psi^0\boxtimes (E\otimes {E'}^*)^{\boxtimes\, d},\;
\pr_{\cY}^!\Psi^0)\; [nd-2d](\frac{nd-2d}{2}),\\
\end{multline*}

\vskip -1.5em\noindent
where the last isomorphism comes from Corolary~\ref{Cor_Hecke}. Since $\pr_{\cY}$ is smooth
of relative dimension $d$, our assertion follows.
\end{Prf}

\smallskip\noindent
\end{Prf}(Proposition~\ref{Pp_really_last})

\smallskip

{\small \scshape Universit\'e Paris-Sud, b\^at. 425, Math\'ematiques, 91405
Orsay France}

e-mail: {\tt Sergey.Lysenko@math.u-psud.fr}

\end{document}